\documentclass[11pt]{amsart}
\usepackage{amssymb}
\usepackage{latexsym}
\usepackage{amsmath}
\usepackage{amsthm}
\usepackage{amsfonts}
\usepackage{color}
\usepackage{pdfsync}
\usepackage[usenames,dvipsnames]{pstricks}
\usepackage{epsfig}
\usepackage{pst-grad} 
\usepackage{pst-plot} 

\newcommand{\R}{\mathbb{R}}
\newcommand{\N}{\mathbb{N}}

\newtheorem{defi}{Definition}[section]

\numberwithin{equation}{section}

\newtheorem{teo}{Theorem}[section]

\newtheorem{proposition}{Proposition}[section]

\newtheorem{lema}{Lemma}[section]
\newtheorem{lemma}{Lemma}[section]

\newtheorem{remark}{Remark}[section]
\newcommand{\bremark}{\begin{remark} \em}
\newcommand{\eremark}{\end{remark} }

\usepackage{amssymb}

\newcommand{\beq}{\begin{equation} }
\newcommand{\eqq}{\end{equation} }
\newcommand{\cuad}{{\sqcap\kern-.68em\sqcup}}

\newcommand{\equ}[1]{(\ref{#1})}

\def\beeq{\begin{equation}}
\def\eeq{\end{equation}}
\newcommand{\begeqaet}{\begin{eqnarray*}}
\newcommand{\eneqaet}{\end{eqnarray*}}

\begin{document}
\begin{center}{\bf  \Large Large solutions to elliptic  equations involving fractional Laplacian  }\medskip

\bigskip


{Huyuan Chen,  Patricio Felmer \\~\\}

 Departamento de Ingenier\'{\i}a  Matem\'atica and
Centro de Modelamiento Matem\'atico
 UMR2071 CNRS-UChile,
 Universidad de Chile\\
 Casilla 170 Correo 3, Santiago, Chile.\\
 {\sl  (hchen@dim.uchile.cl, pfelmer@dim.uchile.cl)} 
  \\~\\and
 
\bigskip

{ Alexander Quaas\\~\\}

Departamento de Matem\'{a}tica, Universidad T\'{e}cnica Federico
Santa Mar\'{i}a\\ Casilla: V-110, Avda. Espa\~{n}a 1680,
Valpara\'{i}so, Chile\\
 {\sl  (alexander.quaas@usm.cl)}

\bigskip

\bigskip

\end{center}

\begin{abstract}

In this article we study existence of boundary blow up solutions for some fractional elliptic equations including 
\begin{eqnarray*}\label{ecbir1a}
(-\Delta)^\alpha u+u^p&=&f\ \
\hbox{in} \ \ \Omega,\\\label{ecbir2a000}
u&=&g\ \ \hbox{on} \ \  \Omega^c,\\\label{ecbir3a000}
\lim_{x\in\Omega, x\to\partial\Omega}u(x)&=&\infty,
\end{eqnarray*}
where $\Omega$ is a bounded domain of class $C^2$, 
$\alpha\in (0,1)$ and the functions  $f:\Omega\to\R$ and $g:
\bar\Omega^c\to\R$ are continuous.... 

We prove existence and uniqueness results fro

ddd

\end{abstract}
\date{}

\setcounter{equation}{0}
\section{ Introduction}
In their pioneering  work,
Keller \cite{K} and Osserman \cite{O} studied the existence of solutions to the nonlinear reaction 
diffusionº equation
\begin{equation}\label{1.1.2}
\left\{ \arraycolsep=1pt
\begin{array}{lll}
 -\Delta u+h(u)=0,\ \ \ \ &
 \mbox{in}\ \ \ &\Omega,\\[2mm]
u=+\infty,&  \mbox{on}\ \ \ &\partial\Omega,
\end{array}
\right.
\end{equation}
where $\Omega$ is a bounded domain in $\R^N$, $N\ge 2$, and $h$ is a nondecreasing positive function. 
They   independently proved that this equation admits a 
solution if and only if $h$ satisfies
\begin{equation}\label{ko}
\int_1^{+\infty}\frac{ds}{\sqrt{H(s)}}<+\infty,
\end{equation}
 where $H(s) = \int_0^{s} h(t)dt$, that in the case of $h(u)=u^p$ means $p>1$. This integral condition on the non-linearity is known as  the  Keller-Osserman
criteria.
The solution of  \equ{1.1.2} found in \cite{K} and \cite{O} exists as  a consequence of  the interaction between the reaction and the difussion term, without the influence of an external source that blows up at the boundary. Solutions exploding at the boundary are usually called boundary blow up solutions or large solutions.  
 From then on, more general boundary blow-up problem:
\begin{equation}\label{3.1.2}
\left\{ \arraycolsep=1pt
\begin{array}{lll}
 -\Delta u(x)+h(x,u)=f(x),\ \ \ \ &
x\in\Omega,\\[2mm]
\lim_{x\in\Omega,\ x\to\partial\Omega}u(x)=+\infty
\end{array}
\right.
\end{equation}
has been extensively studied, see
\cite{AR,BM1,BM2,CCE,DL,DL1,DH,G,LN,MV,MV1,R}. It has being extended
in various ways, weakened the assumptions on the domain and the
nonlinear terms, extended to more general class of equations and
obtained more information on the uniqueness and the asymptotic
behavior of solution at the boundary.

During the last years there has been a renewed and increasing
interest in the study of linear and nonlinear integral operators,
especially, the fractional Laplacian, motivated by great
applications and by important advances on the theory of nonlinear
partial differential equations, see
\cite{CC,CS1,CS2,CF2,FQ0,FQ1,FQ2,FQ3,PVS,S} for details.

In a recent work, Felmer and Quaas \cite{FQ0} considered an analog of \equ{1.1.2} where the laplacian is replaced by the fractional laplacian
\begin{equation}\label{3.1.1}
\left\{ \arraycolsep=1pt
\begin{array}{lll}
 (-\Delta)^{\alpha} u(x)+|u|^{p-1}u=f(x),\ \ \ \ &
\mbox{ in }\quad\Omega,\\[2mm]
u(x)=g(x),& \mbox{ in }\quad \bar\Omega^c,\\[2mm]
\lim_{x\in\Omega,\ x\to\partial\Omega}u(x)=+\infty,
\end{array}
\right.
\end{equation}
where  $\Omega$ is a bounded domain in $\R^N$, $N\ge2$, with  boundary $\partial \Omega$ of class $C^2$, $p>1$ and 
the
fractional Laplacian operator is defined as
$$
(-\Delta)^\alpha u(x)=-\frac12
\int_{\R^N}\frac{\delta(u,x,y)}{|y|^{N+2\alpha}}dy,\ \ x\in \Omega,
$$
with $\alpha\in(0,1)$ and $\delta(u,x,y)=u(x+y)+u(x-y)-2u(x)$.
The authors proved
the existence of a solution  to (\ref{3.1.1}) provided that $g$ 
explodes at the boundary and satisfies other technical conditions. In case the function $g$ blows up with an explosion rate as $d(x)^{\beta}$, with $\beta\in (-\frac{2\alpha}{p-1},0)$ and $d(x)=dist(x,\partial\Omega)$, the solution satisfies
$$
0< \liminf_{x\in\Omega, x\to\partial\Omega}u(x)d(x)^{-\beta}\le \limsup_{x\in\Omega, x\to\partial\Omega}u(x)d(x)^{\frac{2\alpha}{p-1}}<+\infty.
$$
In \cite{FQ0} the explosion is driven by the function $g$. The external source $f$ has a secondary role, not intervening in the explosive character of the solution. $f$ may be bounded or unbounded, in latter case the explosion rate has to be controlled by 
$d(x)^{-2\alpha p/(p-1)}$.

One interesting question not answered in \cite{FQ0} is the existence of a boundary blow up solution without external source, that is assuming $g=0$ in $\bar\Omega^c$ and $f=0$ in $\Omega$, thus extending the original result by
Keller and Osserman, where solutions exists due to the pure interaction between the reaction and the diffusion terms.
It is the purpose of this article  to answer positively this question and to better understand  how the
non-local character influences the large solutions of (\ref{3.1.1})
and what is the structure of the large solutions of (\ref{3.1.1})
with  or without sources. Comparing  with the Laplacian case, where well possedness holds for \equ{3.1.1}, a much richer structure for the solution set appears for the non-local case, depending on the parameters and the data $f$ and $g$. In particular, Theorem \ref{th 3.1.1} shows that existence, uniqueness, non-existence and infinite existence may occur at different values of $p$ and $\alpha$. 

Our first result is on the existence of blowing up solutions driven by the sole interaction between the diffusion and reaction term, assuming the external value $g$ vanishes.
Thus we will be considering the equation
\begin{eqnarray}\label{eteo1}
(-\Delta)^\alpha u+|u|^{p-1}u&=&f\ \
\hbox{in} \ \ \Omega,
\nonumber\\
\label{ecbir2a}
u&=&0\ \ \hbox{in} \ \  \Omega^c,\\
\nonumber
\label{ecbir3a}
\lim_{x\in\Omega, x\to\partial\Omega}u(x)&=&+\infty.
\end{eqnarray}
 On the external source $f$ we will assume the following hypotheses
\begin{itemize}
\item[(H1)] The external source  $f:\Omega\to\R$
is a $C^\beta_{loc}(\Omega)$, for some $\beta>0$.
\item[(H2)] Defining $f_-(x)=\max\{-f(x),0\}$ and $f_+(x)=\max\{f(x),0\}$ we have
\begin{equation*}\label{H21} 
\limsup_{x\in\Omega,x\to\partial\Omega}f_+(x)d(x)^{\frac{2\alpha p}{p-1}}<+\infty \quad\mbox{and}\quad
\lim_{x\in\Omega,x\to\partial\Omega}f_-(x)d(x)^{\frac{2\alpha p}{p-1}}=0.
\end{equation*}
\end{itemize}
A related condition that we need for non-existence results
\begin{itemize}
\item[(H2$^*$)] The function $f$ satisfies
\begin{equation*}\label{H21*} 
\limsup_{x\in\Omega,x\to\partial\Omega}|f(x)|d(x)^{2\alpha}<+\infty.
\end{equation*} 
\end{itemize}
Now we are in a position to state our first theorem
\begin{teo}
\label{th 3.1.1}
Assume that $\Omega$ is an open, bounded and connected domain of
class $C^2$ and $\alpha\in(0,1)$. Then we have:

\noindent
{\bf Existence}: Assume that $f$ satisfies (H1) and (H2), then there exists $\tau_0(\alpha)\in (-1,0)$ such that for every $p$ satisfying
\begin{equation}\label{pp1}
1+2\alpha<p<1-\frac{2\alpha}{\tau_0(\alpha)},
\end{equation}
the equation \equ{eteo1} possesses at least one solution 
 $u$ satisfying
\begin{equation}\label{3.1.5}
0<\liminf_{x\in\Omega,x\to\partial\Omega}
u(x)d(x)^{\frac{2\alpha}{p-1}}\le\limsup_{x\in\Omega,x\to\partial\Omega}
u(x)d(x)^{\frac{2\alpha}{p-1}} <+\infty.
\end{equation}
{\bf Uniqueness}: If $f$ further satisfies $f\ge 0$ in $\Omega$, then $u>0$ in $\Omega$ and $u$ is the unique
solution of \equ{eteo1} satisfying \equ{3.1.5}.

\noindent
{\bf Nonexistence}: If $f$ satisfies (H1) and (H2$^*$), then in the following three cases:
\begin{itemize}
\item[i)] For any $\tau\in(-1,0)\setminus\{-\frac{2\alpha}{p-1},\
\tau_0(\alpha)\}$ and $p$ satisfying \equ{pp1} or
\item[ii)] For any $\tau\in(-1,0)$ and 
\begin{equation}\label{pp2}
p\ge 1-\frac{2\alpha}{\tau_0(\alpha)}\mbox{ or}
\end{equation}
\item[iii)] For any $\tau\in(-1,0)\setminus\{\tau_0(\alpha)\}$ and 
\begin{equation}\label{pmalpha}
1<p\le 1+2\alpha,
\end{equation}
\end{itemize}
equation  \equ{eteo1} does not have a solution $u$ satisfying
\begin{equation}\label{1.7}
0<\liminf_{x\in\Omega,x\to\partial\Omega}
u(x)d(x)^{-\tau}\le\limsup_{x\in\Omega,x\to\partial\Omega}
u(x)d(x)^{-\tau} <+\infty.
\end{equation} 
\noindent
{\bf Special existence for $\tau=\tau_0(\alpha)$.}
 Assume $f(x)\equiv0,\ x\in\Omega$ and
that
\begin{equation}\label{observartion}
\max\{1-\frac{2\alpha}{\tau_0(\alpha)}+\frac{\tau_0(\alpha)+1}{\tau_0(\alpha)},1\}<p<1-\frac{2\alpha}{\tau_0(\alpha)}.
\end{equation}
Then, there exist constants $C_1\geq 0$ and $C_2>0$, such that for any $t>0$ there is  a  positive solution $u$ of equation
(\ref{eteo1})  satisfying
\begin{equation}\label{obs3.1.5}
C_1d(x)^{\min\{\tau_0(\alpha)p+2\alpha,0\}}\leq td(x)^{\tau_0(\alpha)}-u(x)\le
C_2 d(x)^{\min\{\tau_0(\alpha)p+2\alpha,0\}}.
\end{equation}
\end{teo}


\begin{remark}
We remark that hypothesis (H2) and ($\rm{H2^*}$) are satisfied when $f\equiv 0$, so this theorem answer the question on existence rised in \cite{FQ0}.
We also observe that a function $f$ satisfying (H2) may also satisfy
$$\lim_{x\in\Omega,x\in\partial\Omega}f(x)=-\infty,$$
what matters is that the rate of explosion is smaller than $\frac{2\alpha p}{p-1}$.
\end{remark}

For proving the existence part of this theorem we will construct appropriate super and sub-solutions. This construction involves the one dimensional truncated laplacian of power functions given by 
\begin{equation}\label{3.1.4}
C(\tau)=\int^{+\infty}_{0}\frac{\chi_{(0,1)}(t)|1-t|^{\tau}+(1+t)^{\tau}-2}{t^{1+2\alpha}}dt,
\end{equation}
for $\tau\in (-1,0)$ and where $\chi_{(0,1)}$ is the characteristic function of the interval $(0,1)$. The number  $\tau_0(\alpha)$ appearing in the statement of our theorems is precisely the unique 
$\tau\in (-1,0)$ satisfying $C(\tau)=0$. See Proposition \ref{pr 3.3.2} for details.

\begin{remark}
For the uniqueness, we would like to
mention that, by using iteration technique, Kim in \cite{K1} has
proved the uniqueness of solution to the problem
\begin{equation}\label{1.1.1}
\left\{ \arraycolsep=1pt
\begin{array}{lll}
 -\Delta u+u_+^p=0,\ \ \ \ &
 \mbox{in}\ \ \ &\Omega,\\[2mm]
u=+\infty,&  \mbox{in}\ \ \ &\partial\Omega,
\end{array}
\right.
\end{equation}
where $u_+=\max\{u,0\}$, under the hypotheses that $p>1$ and
$\Omega$ is bounded and satisfying
$\partial\Omega=\partial\bar\Omega$.  Garc\'{i}a-Meli\'{a}n in
\cite{G,GG} introduced some improved  iteration technique to obtain
the uniqueness for problem (\ref{1.1.1}) with replacing nonlinear
term by $a(x)u^p$. However, there is a big difficulty for us to
extend the iteration technique to our problem (\ref{3.1.1})
involving fractional Laplacian, which  is caused  by the nonlocal
character.
\end{remark}

\medskip

In the second part, we are also interested in considering the existence of blowing up solutions driven by external source $f$ on which we assume 
 the following hypothesis
\begin{itemize}
\item[(H3)]
There exists $\gamma\in(-1-2\alpha,0)$ such that
$$
0<\liminf_{x\in\Omega,x\to\partial\Omega}f(x)d(x)^{-\gamma}\leq
\limsup_{x\in\Omega,x\to\partial\Omega}f(x)d(x)^{-\gamma}<+\infty.
$$
\end{itemize}
Depending on the size of $\gamma$ we will say that the external source is weak or strong. 
In order to gain in clarity, in this case we will state separately the existence, uniqueness and non-existence theorem in this source-driven case. 
\begin{teo} [Existence]
\label{th 3.1.2}
Assume that $\Omega$ is an open, bounded and connected domain of
class $C^2$. Assume that $f$ satisfies (H1)  and let $\alpha\in(0,1)$ then we have:

$(i)$ (weak source) If $f$ satisfies (H3) with 
\begin{equation}\label{gamma1}
-2\alpha-\frac{2\alpha}{p-1}\leq\gamma<-2\alpha,
\end{equation}
then, for every $p$ such that \equ{pp2} holds,
equation  \equ{ecbir2a} possesses at
least one solution $u$,  with asymptotic behavior near the boundary given by
\begin{equation}\label{3.1.10}
0<\liminf_{x\in\Omega,x\to\partial\Omega}
u(x)d(x)^{-\gamma-2\alpha}\le\limsup_{x\in\Omega,x\to\partial\Omega}
u(x)d(x)^{-\gamma-2\alpha} <+\infty.
\end{equation}
$(ii)$ (strong source)  If $f$ satisfies (H3) with 
\begin{equation}\label{gamma2}
-1-2\alpha<\gamma<-2\alpha-\frac{2\alpha}{p-1}
\end{equation}
then, for every $p$ such that
\begin{equation}\label{pp4} 
p>1+2\alpha,
 \end{equation}
 equation  \equ{ecbir2a} possesses at
least one solution $u$,  with asymptotic behavior near the boundary given by 
\begin{equation}\label{3.1.11}
0<\liminf_{x\in\Omega,x\to\partial\Omega}
u(x)d(x)^{-\frac{\gamma}{p}}\le\limsup_{x\in\Omega,x\to\partial\Omega}
u(x)d(x)^{-\frac{\gamma}{p}} <+\infty.
\end{equation}
\end{teo}

As we already mentioned, in Theorem \ref{th 3.1.1} the existence of blowing up solutions results from the interaction between the reaction $u^p$ and the diffusion term $(-\Delta)^\alpha$, while the role of the external source $f$ is secondary. In contrast, in  Theorem \ref{th 3.1.2} the existence of blowing up solutions results on the interaction between the external source, and the diffusion term in case of weak source and the interaction between the  external source and the reaction term in case of strong source.

Regarding uniqueness result for solutions of \equ{ecbir2a}, as in Theorem  \ref{th 3.1.1} we will assume that $f$ is non-negative, hypothesis that we need for  technical reasons.  We have

\begin{teo} [Uniqueness]
\label{th 1.1}
Assume that $\Omega$ is an open, bounded and connected domain of
class $C^2$, $\alpha\in(0,1)$ and
   $f$ satisfies (H1) and $f\ge 0$. Then  we have
\begin{itemize}
\item[i)] (weak source)  the solution of \equ{ecbir2a} satisfying  (\ref{3.1.10}) is positive and unique, and 
\item[ii)] (strong source)  the solution of \equ{ecbir2a} satisfying (\ref{3.1.11}) is positive and unique.
\end{itemize}
\end{teo}

We complete our theorems with a non-existence result for solution with a previously defined asymptotic behavior, as we saw in Theorem \ref{th 3.1.1}.
We have
\begin{teo}[Non-existence]
\label{nonexistence}
Assume that $\Omega$ is an open, bounded and connected domain of
class $C^2$, $\alpha\in(0,1)$ and $f$ satisfies
$(H1)$, $(H3)$ and $f\ge 0$. Then we have

\begin{itemize}
\item[i)] (weak source) Suppose that $p$ satisfies \equ{pp2}, $\gamma$ satisfies \equ{gamma1} and $\tau\in(-1,0)\setminus\{\gamma+2\alpha\}$. Then 
equation \equ{eteo1} does not have a solution  $u$ satisfying \equ{1.7}.  
%

\item[ii)]  (strong source) Suppose that $p$ satisfies \equ{pp4}, $\gamma$ satisfies \equ{gamma2} and $\tau\in(-1,0)\setminus\{\frac\gamma p\}$. Then, 
equation \equ{eteo1} does not have a solution  $u$ satisfying \equ{1.7}.  
\end{itemize}
\end{teo}

 All theorems stated so far deal with equation \equ{3.1.1} in the case $g\equiv 0$, but they may also be applied when  $g\not\equiv 0$ and, in particular, these result improve those given in \cite{FQ0}. In what follows we describe how to obtain this. 
We start with some notation, we consider $
L^1_\omega(\bar\Omega^c)$ the weighted $L^1$ space in $\bar\Omega^c$ with weight
$$
\omega(y)= \frac{1}{1+|y|^{N+2\alpha}},\quad\mbox{for all } y\in\R^N.
$$  
Our hypothesis  on the  external values  $g$ is the following
\begin{itemize}
\item[$(H4)\ $]
\begin{enumerate}\item[]
The function $g:\bar\Omega^c\to\R$ is measurable  and $g\in
L^1_\omega(\bar\Omega^c)$.
\end{enumerate}
\end{itemize}
Given $g$ satisfying $(H4)$, we define
\begin{equation}\label{G} G(x)=\frac12 \int_{\R^N}\frac{\tilde
g(x+y)}{|y|^{N+2\alpha}}dy,\ \ x\in\Omega,\end{equation} where
\begin{equation}\label{g}
\tilde g(x)=\left\{ \arraycolsep=1pt
\begin{array}{lll}
 0,\ \ \ \ &
x\in \bar\Omega,\\[2mm]
 g(x),\ &x\in\bar\Omega^c.
\end{array}
\right.
\end{equation}
We observe that 
$$
G(x)=-(-\Delta)^{\alpha} \tilde g(x),\ \ x\in
\Omega.
$$
 Hypothesis $(H4)$ implies that $G$ is
continuous in $\Omega$ as seen in Lemma \ref{lm 3.2.1} and has an explosion of order $d(x)^{\beta-2\alpha}$
towards the boundary $\partial \Omega$, if $g$ has an explosion of order
$d(x)^{\beta}$ for some $\beta\in (-1,0)$, as we shall see in Proposition \ref{lm 3.3.21}.
We observe that under  the hypothesis $(H4)$, if $u$ is a  solution of equation (\ref{3.1.1}), then
$u-\tilde g$ is the solution of
\begin{equation}\label{without outside source}
\left\{ \arraycolsep=1pt
\begin{array}{lll}
 (-\Delta)^{\alpha} u(x)+|u|^{p-1}u(x)=f(x)+G(x),\ \ \ \ &
x\in\Omega,\\[2mm]
u(x)=0,& x\in\bar\Omega^c,\\[2mm]
\lim_{x\in\Omega,\ x\to\partial\Omega}u(x)=+\infty
\end{array}
\right.
\end{equation}
and vice versa,  if $v$ is a solution of (\ref{without outside
source}), then $v+\tilde g$ is a solution of (\ref{3.1.1}).

Thus, using Theorem  \ref{th 3.1.1}-\ref{nonexistence}, we can state the corresponding results of existence, uniqueness and non-existence for (\ref{3.1.1}), combining  $f$
with $g$ to define a new external source
\begin{equation}\label{F}
 F(x)= G(x)+f(x),\ \ \ x\in \Omega.
\end{equation}
With this we can state appropriate hypothesis for $g$ and thus we can write theorems, one  corresponding to each of Theorem 1.1, 1.2, 1.3 and 1.4.
Even though, at first sight we need that $G(x)$ is $C^\beta_{loc}(\Omega)$, actually continuity of $g$ is sufficient, as we discuss  Remark \ref{RQQQ}.

Moreover, in  Remark \ref{RQQQ1} we explain how   our results in this paper allows to give a different proof of those  obtained by Felmer and Quaas in \cite{FQ0}, generalizing them.

\medskip

 This article is organized as
follows. In Section \S2 we present some preliminaries to introduce
the notion of viscosity solutions, comparison and stability theorems
in case of explosion at the boundary. Then we prove an existence theorem for the nonlinear problem with blow up at the boundary, assuming the existence of ordered. Section \S3 is devoted to obtain crucial estimates
used to construct super and sub-solutions.  In Section \S4
 we prove the existence of solution to (\ref{eteo1}) in Theorem \ref{th 3.1.1} and Theorem \ref{th 3.1.2}.
 In section \S5, we give the proof
of the uniqueness of solution to (\ref{eteo1}) in Theorem \ref{th
3.1.1} and Theorem \ref{th 1.1}. Finally, the nonexistence related
to Theorem \ref{th 3.1.1} and  Theorem \ref{nonexistence} are shown in Section \S6.

\setcounter{equation}{0}
\section{Preliminaries and existence theorem}

The purpose of this section is to introduce some preliminaries and prove  an  existence  theorem for blow-up solutions assuming  the existence of ordered   super-solution and sub-solution which  blow up at the boundary. In order to prove this theorem we adapt the  theory of viscosity to allow for boundary blow up. 

 We start this section by defining   the notion of viscosity solution for non-local equation, allowing  blow up at the boundary, see for example \cite{CS2}. We consider the equation of the form:
 \begin{equation}\label{eqhh}
(-\Delta)^{\alpha}  u=h(x, u) \quad\mbox{in} \quad\Omega,\quad
u=g \quad\mbox{in} \quad\Omega^c.
\end{equation}
\begin{defi} \label{de 2.2}
We say that a function $u:(\partial \Omega)^c\to \R$, continuous   in $\Omega$ and in $
L^1_\omega(\R^N)$ is a viscosity super-solution
(sub-solution) of (\ref{eqhh}) if
$$u\geq g\ (\mbox{resp.}\ u\leq g)\ \mbox{in}\ \bar \Omega^c$$
and for every point $x_0\in\Omega$ and some  neighborhood $V$ of
$x_0$ with $\bar V\subset \Omega $ and for any $\phi\in C^2(\bar V)$
such that $u(x_0)=\phi(x_0)$ and
$$u(x)\ge \phi(x)\ (\mbox{resp.}\ u(x)\leq \phi(x))\ \mbox{for\ all}\ x\in V,
$$
 defining
\begin{eqnarray*}
\tilde u=\left\{ \arraycolsep=1pt
\begin{array}{lll}
\phi\ \ & \mbox{in}\ \ &V,\\[2mm]
u & \mbox{in}\ \ &V^c,
\end{array}
\right.
\end{eqnarray*}
we have
$$(-\Delta)^{\alpha} \tilde u(x_0)\geq h(x_0, u(x_0))\
(\mbox{resp.}(-\Delta)^{\alpha} \tilde u(x_0)\leq h(x_0, u(x_0)).$$
We say that $u$ is a  viscosity solution of (\ref{eqhh}) if it is a viscosity super-solution and also a  viscosity sub-solution  of (\ref{eqhh}).
\end{defi}
It will be convenient for us to have also a notion of classical solution.
\begin{defi} \label{de 2.3}
We say that a function $u:(\partial \Omega)^c\to \R$, continuous   in $\Omega$ and in $
L^1_\omega(\R^N)$ is a classical solution of \equ{eqhh} if $(-\Delta)^\alpha u(x)$ is well defined for all $x\in\Omega$,
$$
(-\Delta)^\alpha u(x)=h(x,u(x)),\quad\mbox{for all } x\in\Omega
$$
and $u(x)=g(x)$ a.e. in $\overline{\Omega}^c$. Classical super and sub-solutions are defined similarly.
\end{defi}
\noindent
Next we have our first regularity theorem.
\begin{teo}\label{regularity}
 Let $g\in L^1_\omega(\R^N)$ and $f\in C^\beta_{loc}(\Omega)$, with $\beta\in (0,1)$, and $u$ be a viscosity solution of
$$
(-\Delta^\alpha u)=f \quad\mbox{ in }\quad \Omega,\quad
u=g \quad\mbox{in} \quad\Omega^c,
$$
 then there exists $\gamma>0$ such that $u\in C^{2 \alpha+\gamma}_{loc}(\Omega)$  
\end{teo}
\noindent
{\bf Proof.}  Suppose without loss of generality that $ B_1\subset \Omega$  and $f \in C^\beta (B_1)$.
 Let  $\eta$ be a non-negative, smooth function with support in $B_1$, such that $\eta=1$ in $B_{1/2}$.
Now we look at the equation
$$
-\Delta w=\eta f \quad\mbox{ in }\quad \ \R^N.
$$
 By H\"older regularity theory for the Laplacian we find $w\in C^{2,\beta}$,
so that  $(-\Delta)^{1-\alpha}w\in  C^{2\alpha+\beta}$, see \cite{STEIN} or  Theorem 3.1 in \cite{FQT}.
Then, since
$$
(-\Delta)^\alpha (u-(-\Delta)^{1-\alpha}w)=0 \quad\mbox{ in }\quad B_{1/2} ,
$$
we can use Theorem 1.1 and Remark 9.4 of \cite{CS_Evans}  (see also Theorem 4.1 there), to obtain that there exist $\tilde\beta$ such that
$u-(-\Delta)^{1-\alpha}w \in C^{2\alpha+\tilde\beta}(B_{1/2})$, from where we conclude.\hfill$\Box$

\medskip

The Maximum and the Comparison Principles are key tools in the analysis, we present them here for completitude.
\begin{teo}\label{MP}(Maximum principle)
Let ${\mathcal{O}}$ be an open and  bounded  domain of $\R^N$ and $u$
be a classical solution of
\begin{equation}\label{2.2}
(-\Delta)^\alpha u\le 0 \ \ \ \mbox{in}\ \ \ {\mathcal{O}},
\end{equation}
continuous in $\bar{{\mathcal{O}}}$ and bounded from above in $\R^N$.
Then $u(x)\le M,$ for all  $x\in {\mathcal{O}},$ where
$M=\sup_{x\in{\mathcal{O}}^c}u(x)<+\infty.$
\end{teo}
\noindent
{\bf Proof.} If the conclusion is false, then there exists $x'\in
{\mathcal{O}}$ such that $u(x')>M$. By continuity of $u$, there exists
$x_0\in {\mathcal{O}}$ such that
$$u(x_0)=\max_{x\in{\mathcal{O}}}u(x)=\max_{x\in\R^N}u(x)$$ 
and then  $(-\Delta)^\alpha u(x_0)>0 $, which contradicts  (\ref{2.2}).
\hfill$\Box$

\medskip

\begin{teo}(Comparison Principle)\label{comparison}
Let $u$ and $v$ be classical super-solution and
sub-solution of
$$
 (-\Delta)^{\alpha} u+h(u)=f\ \ \mbox{in}\ \ {\mathcal{O}},$$ 
respectively, where ${\mathcal{O}}$ is an open, bounded domain,
 the functions $f:{\mathcal{O}}\to\R$ is
continuous and $h:\R\to\R$ is increasing.
Suppose further that $u$ and $v$ are
continuous in $\bar {\mathcal{O}}$ and $v(x)\le u(x)$ for all  $x\in{\mathcal{O}}^c$. Then
$$u(x)\ge v(x),\ x\in {\mathcal{O}}.$$
\end{teo}
\noindent
{\bf Proof.} Suppose by contradiction that $w=u-v$ has a negative minimum in $x_0\in{\mathcal{O}}$, then  $(-\Delta)^{\alpha} w(x_0)<0$ and so, by assumptions on $u$ and $v$,  $h(u(x_0))>h(v(x_0))$, which contradicts the monotonicity of $h$.  \hfill$\Box$\\

We devote the rest of the section to the proof of the existence theorem through super and sub-solutions. We prove the theorem by an approximation procedure for which we need some preliminary steps. We need to deal with a Dirichlet problem  involving fractional laplacian operator and with
exterior data which blows up away from the boundary. Precisely, on the exterior data $g$, we assume the following hypothesis, given an open, bounded set  ${\mathcal{O}}$ in $\R^N$ with $C^2$ boundary:
\begin{itemize}
\item[$(G)\ $] $g:{\mathcal{O}^c} \to \R$ is in $L^1_\omega({\mathcal{O}}^c)$  and it is of class  $C^2$ in $\{z\in{\mathcal{O}}^c,
dist(z,\partial{\mathcal{O}})\leq \delta \}$,   where
$\delta>0$.
\end{itemize}

In studying  the nonlocal problem (\ref{3.1.1}) with 
explosive exterior source, we have to adapt  the stability theorem and the existence theorem for the linear Dirichlet
problem. The following lemma is important in this direction. 
\begin{lema}\label{lm 3.2.1}
Assume that ${\mathcal{O}}$ is an open, bounded domain in $\R^N$ with 
 $C^2$ boundary. Let $w:\R^N\to\R$:
 \\ $(i)$ If 
$w\in L^1_\omega(\R^N)$ and $w$ is of class $C^2$ in $\{z\in\R^N,
d(z, \mathcal{O})\leq \delta\}$ for some $\delta>0$, then $(-\Delta)^{\alpha} w$
is  continuous in $\bar{{\mathcal{O}}}$.\\
$(ii)$ If  $w\in L^1_\omega(\R^N)$ and $w$ is of class $C^2$ in ${\mathcal{O}}$,
then $(-\Delta)^{\alpha} w$ is  continuous in ${\mathcal{O}}$.
\\
$(iii)$ If  $w\in L^1_\omega(\R^N)$ and $w\equiv 0$ in ${\mathcal{O}}$, then $(-\Delta)^{\alpha} w$ is  continuous in ${\mathcal{O}}$.
\end{lema}
\noindent
{\bf Proof.}  We first prove (ii). Let $x\in\Omega$ and $\eta>0$ such that $B(x,2\eta)\subset \Omega$. Then we consider 
$$
(-\Delta)^\alpha u(x)=L_1(x)+L_2(x),
$$
where
$$
L_1(x)=\int_{B(0,\eta)}\frac{\delta(u,x,y)}{|y|^{N+2\alpha}}dy \quad\mbox{and}\quad 
L_2(x)=\int_{B(0,\eta)^c}\frac{\delta(u,x,y)}{|y|^{N+2\alpha}}dy.
$$
Since $w$ is of class $C^2$ in $\mathcal{O}$, we may write $L_1$ as
$$
L_1(x)=\int_0^\eta\left\{ \int_{S^{N-1}}\int_{-1}^1\int_1^1t\omega^tD^2w(x+str\omega)\omega dtdsd\omega\right\} r^{1-\alpha} dr,
$$
where the term inside the brackets is  uniformly continuous in $(x,r)$, so  the resulting function $L_1$ is continuous. On the other hand we may write $L_2$ as
$$
L_2(x)=-2w(x)\int_{B(0,\eta)^c}\frac{dy}{|y|^{N+2\alpha}}-2\int_{B(x,\eta)^c}\frac{w(z)dz}{|z-x|^{N+2\alpha}},
$$
from where $L_2$ is also continuous.
The proof of (i) and (iii) are similar. 
\hfill$\Box$

\medskip

The next theorem gives the stability property for viscosity solutions in our setting. 
\begin{teo}\label{stability}
Suppose that ${\mathcal{O}}$ is an open, bounded and $C^2$ domain and
$h:\R\to \R$ is continuous.  Assume that $(u_n)$, $n\in\N$ is a sequence of functions,   bounded in
$L^1_\omega({\mathcal{O}}^c)$ and $f_n$ and $f$
are continuous in ${\mathcal{O}}$ such that:

$(-\Delta)^\alpha u_n+ h(u_n)\ge f_n\ (\mbox{resp.}\ (-\Delta)^\alpha
u_n+ h(u_n)\le f_n)$ in ${\mathcal{O}}$ in viscosity sense,

$u_n\to u$ locally uniformly in  ${\mathcal{O}}$,

$ u_n\to u$  in  $L^1_\omega(\R^N)$, and

$ f_n\to f$ locally uniformly in  ${\mathcal{O}}$.\\
    Then, $(-\Delta)^\alpha u+h(u)\ge f\ (\mbox{resp.}\ (-\Delta)^\alpha u+h(u)\le f)$ in ${\mathcal{O}}$ in viscosity sense.
\end{teo}
\noindent
{\bf Proof.} If $|u_n|\leq C$ in ${\mathcal{O}}$ then we use Lemma 4.3 of \cite{CS2}.
If $u_n$  is unbounded in ${\mathcal{O}}$, then $u_n$ is bounded in 
${\mathcal{O}}_k=\{x\in {\mathcal{O}}, dist(x,\partial{\mathcal{O}})\ge\frac1k\}$, since $u_n$ is continuous in ${\mathcal{O}}$, and then by  Lemma 4.3 of \cite{CS2}, $u$
is a viscosity solution of $(-\Delta)^\alpha u+h(u)\ge f$ in
${\mathcal{O}}_k$ for any $k$. Thus $u$ is a viscosity solution of
$(-\Delta)^\alpha u+h(u)\ge f$ in ${\mathcal{O}}$ and the proof is completed.
\hfill$\Box$\\

An existence result for the Dirichlet problem is given as follows:

\begin{teo}\label{th 3.2.3}
 Suppose that ${\mathcal{O}}$ is an open, bounded and $C^2$ domain,  $g:{\mathcal{O}}^c\to\R$ satisfies   $(G)$, $f:\bar{{\mathcal{O}}}\to\R$ is continuous, $f\in C^\beta_{loc}({\mathcal{O}})$, with $\beta\in(0,1)$,  and $p>1$.
Then there exists a classical solution $u$ of
\begin{equation}\label{3.2.4}
\left\{ \arraycolsep=1pt
\begin{array}{lll}
 (-\Delta)^{\alpha} u(x)+|u|^{p-1}u(x)=f(x),\ \ \ \ &
x\in{\mathcal{O}},\\[2mm]
u(x)=g(x),\ &x\in{\mathcal{O}}^c,
\end{array}
\right.
\end{equation}
which is continuous in $\bar{{\mathcal{O}}}$.
\end{teo}
In proving Theorem \ref{th 3.2.3}, we will use the following
lemma:
\begin{lema}\label{lm 3.2.2}
Suppose that ${\mathcal{O}}$ is an open, bounded and $C^2$ domain, $f:\bar{{\mathcal{O}}}\to\R$ is continuous  and $C>0$. Then there
exist a classical solution of
\begin{eqnarray}
\left\{ \arraycolsep=1pt
\begin{array}{lll}\label{xxx}
 (-\Delta)^{\alpha} u(x)+Cu(x)=f(x),\ \ \ \ &
x\in{\mathcal{O}},\\[2mm]
u(x)=0,\ &x\in{\mathcal{O}}^c,
\end{array}
\right.
\end{eqnarray}
which is continuous in $\bar{{\mathcal{O}}}$.
\end{lema}
\noindent
{\bf Proof.} For the existence of a viscosity solution $u$ of \equ{xxx}, that is continuous in $\bar{{\mathcal{O}}}$, we refers to  Theorem 3.1 in
\cite{FQ0}.   Now we apply Theorem 2.6 of \cite{CS2} to conclude that $u$ is $C^\theta_{loc}({\mathcal{O}})$, with $\theta>0$, and then we use Theorem \ref{regularity} to conclude that the solution is classical (see also Proposition 1.1 and 1.4 in \cite{rosoton}).   \hfill$\Box$

\medskip

Using Lemma \ref{lm 3.2.2}, we find  $\bar V$,  a classical solution of
\begin{equation} \label{3.2.V}
\left\{ \arraycolsep=1pt
\begin{array}{lll}
 (-\Delta)^{\alpha} \bar V(x)=-1,\ \ \ \ &
x\in{\mathcal{O}},\\[2mm]
\bar V(x)=0,\ &x\in{\mathcal{O}}^c,
\end{array}
\right.
\end{equation}
which is continuous in $\bar{{\mathcal{O}}}$ and  negative in
${\mathcal{O}}$. it is classical since we apply Theorem 2.6 of \cite{CS2} to conclude that $u$ is $C^\theta_{loc}({\mathcal{O}})$, with $\theta>0$, and then we use Theorem \ref{regularity} to conclude that the solution is classical (see also Proposition 1.1 and 1.4 in \cite{rosoton}).

Now we prove Theorem \ref{th 3.2.3}.

\noindent 
{\bf Proof of Theorem \ref{th 3.2.3}.}
Under  assumption $(G)$ and in view of the hypothesis on $\mathcal{O}$, we may extend $g$ to $\bar g$ in
$\R^N$ as a 
$C^2$  function in $\{z\in \R^N,d(z,\mathcal{O})\leq \delta\}$.   We certainly have $\bar g\in
L^1_\omega(\R^N)$ and, by Lemma \ref{lm 3.2.1}
 $(-\Delta)^{\alpha} \bar g$ is continuous in
$\bar{{\mathcal{O}}}$. 
Next we use 
  Lemma \ref{lm 3.2.2} to find a solution $v$ of equation \equ{xxx} with $f(x)$ replaced by $f(x)-(-\Delta)^\alpha \bar g(x)-C\bar g(x)$, where $C>0$. 
Then
we define $u=v+\bar g$ and we see that $u$ is continuous in  $\bar{\mathcal{O}}$ and it satisfies in the viscosity sense  
\begin{eqnarray*} \left\{
\arraycolsep=1pt
\begin{array}{lll}
 (-\Delta)^{\alpha} u(x)+Cu(x)=f(x),\ \ \ \ &
x\in{\mathcal{O}},\\[2mm]
u(x)=g(x),\ &x\in{\mathcal{O}}^c.
\end{array}
\right.
\end{eqnarray*}
Now we use Theorem Theorem 2.6 in \cite{CS2} and then Theorem \ref{regularity}  to conclude that  $u$ is a classical solution. 
Continuing the proof, we find super and sub-solutions for \equ{3.2.4}.
We define 
$$ u_\lambda(x)=\lambda \bar V(x)+\bar g(x), \
x\in\R^N,
$$ 
where $\lambda\in\R$ and $\bar V$ is given in 
(\ref{3.2.V}). We see that $ u_\lambda(x)=g(x)$  in
${\mathcal{O}}^c$ for any $\lambda$ and for $\lambda$ large (negative), $u_\lambda$ satisfies
\begin{eqnarray*}
(-\Delta)^{\alpha} u_\lambda(x)+
|u_\lambda(x)|^{p-1}u_\lambda(x)-f(x)\ge(-\Delta)^{\alpha}\bar g(x)-\lambda-f(x)-|\bar g(x)|^p,
\end{eqnarray*} 
for $x\in {\mathcal{O}}$. Since 
$(-\Delta)^{\alpha}\bar g$, $\bar g$ and $f$ are
bounded in $\bar {\mathcal{O}}$, choosing $\lambda_1<0$ large enough we find 
 that $u_{\lambda_1}\ge0$ is a super-solution of (\ref{3.2.4})
with $ u_{\lambda_1}=g$ in ${\mathcal{O}}^c$.

On the other hand, for $\lambda>0$ we have
\begin{eqnarray*}
(-\Delta)^{\alpha} u_\lambda(x)+
|u_\lambda|^{p-1}u_\lambda(x)-f(x)\le(-\Delta)^{\alpha}\bar
g(x)-\lambda+|\bar g|^{p-1}\bar g(x) -f(x).
\end{eqnarray*}
As before, there is $\lambda_2>0$ large enough, so that  $u_{\lambda_2}$ is a sub-solution of (\ref{3.2.4})
with $ u_{\lambda_2}=g$ in ${\mathcal{O}}^c$. Moreover, we have that $u_{\lambda_2}< u_{\lambda_1}$ in ${\mathcal{O}}$ and
$u_{\lambda_2}=u_{\lambda_1}=g$ in ${\mathcal{O}}^c$.

Let $u_0=u_{\lambda_2}$ and  define iteratively, using the above argument, the
sequence of functions $u_n\ (n\ge1)$ as the classical solutions of
\begin{equation*}\label{3.2.6} 
\arraycolsep=1pt
\begin{array}{lll}
 (-\Delta)^{\alpha} u_n(x)+Cu_n(x)=f(x)+Cu_{n-1}(x)-|u_{n-1}|^{p-1}u_{n-1}(x),\  &
x\in{\mathcal{O}},\\ 
\hspace{2.8cm}u_n(x)=g(x),\quad   x\in{\mathcal{O}}^c,   &\hspace{2.8cm}~
\end{array}
\end{equation*}
where $C>0$ is  so that the function
$r(t)=Ct-|t|^{p-1}t$ is increasing in the interval $[\min_{x\in\bar
{\mathcal{O}}} u_{\lambda_2}(x),\max_{x\in\bar
{\mathcal{O}}}u_{\lambda_1}(x)]$. 
Next, using Theorem \ref{comparison} we get 
$$u_{\lambda_2}\leq u_n\leq u_{n+1}\leq u_{\lambda_1}\ \ \mbox{in}\ {\mathcal{O}}, \quad \mbox{for all } n\in \N.$$
Then we define 
$u(x)=\lim_{n\to+\infty}u_n(x),$ for $x\in{\mathcal{O}}$ and $
u(x)=g(x),$ for $x\in {\mathcal{O}}^c$ and we have
\begin{equation}\label{banda}
u_{\lambda_2}\leq u\leq u_{\lambda_1}\ \ \mbox{in}\ \ {\mathcal{O}}.
\end{equation}
Moreover, $u_{\lambda_1}, u_{\lambda_2}\in L^1_\omega(\R^N)$
so that  $
u_n\to u$ in $L^1_\omega(\R^N),$ as $n\to\infty$.

By interior estimates as given in \cite{CS1},  for any
compact set $K$ of ${\mathcal{O}}$, we have that $u_n$ has uniformly bounded $C^\theta(K)$ norm. So,  by Ascoli-Arzela Theorem we have that $u$
is continuous in $K$ and $u_n\to u$ uniformly in $K$. Taking  a
sequence of compact sets $K_n=\{z\in{\mathcal{O}}, d(z,\partial {\mathcal{O}})\ge
\frac1n\}$, and ${\mathcal{O}}=\cup^{+\infty}_{n=1}K_n,$  we find that $u$ is
continuous in ${\mathcal{O}}$ and, by Theorem \ref{stability}, $u$ is a viscosity solution of
(\ref{3.2.4}). Now we apply Theorem 2.6 of \cite{CS2} to find that u is $C^\theta_{loc}({\mathcal{O}})$, and then we  use Theorem \ref{regularity} con conclude that $u$ is a classical  solution.  
In addition,  $u$ is continuous up to the boundary by \equ{banda}. 
\hfill $\Box$

Now we are in a position to prove the main theorem of this section. We prove the existence of a blow-up solution of \equ{eteo1} assuming the existence of  suitable super and sub-solutions.
\begin{teo}\label{th 3.2.4}
  Assume that $\Omega$ is an open, bounded domain of
class $C^2$, $p>1$ and $f$ satisfy $(H1)$.   Suppose there exists
a super-solution $\bar U$ and a sub-solution $\underline U$ of
(\ref{eteo1}) such that $\bar U$ and $\underline U$ are of class $C^2$
 in $\Omega$, $\underline U$, $\bar U\in 
L^1_\omega(\R^N)$,
$$\bar U\geq \underline U\ \ \mbox{in}\ \Omega,\ \
\liminf_{x\in\Omega,x\to\partial\Omega}\underline U(x)=+\infty\ \
\mbox{and}\ \ \bar U= \underline U=0\ \ \mbox{in}\ \bar\Omega^c.
$$ 
Then
there exists at least one solution $u$ of (\ref{eteo1}) in the
viscosity sense and
$$\underline U\leq u\leq\bar U \ \ \mbox{in}\ \ \Omega.
$$
Additionally, if $f\ge0$ in $\Omega,$ then
$u>0$ in $\Omega$.
\end{teo}
\noindent
{\bf Proof.} Let us consider $\Omega_n=\{x\in\Omega: d(x)>1/n\}$ and use Theorem \ref{th 3.2.3} to find a solution $u_n$ of
\begin{equation}\label{3.2.7} \left\{ \arraycolsep=1pt
\begin{array}{ll}
 (-\Delta)^{\alpha} u(x)+|u|^{p-1}u(x)=f(x),\ \ \ \ &
x\in\Omega_{n},\\[2mm]
u(x)=\underline U(x),\ &x\in \Omega_{n}^c,
\end{array}
\right.
\end{equation}
 We just replace
${\mathcal{O}}$ by $\Omega_{n}$ and define
$\delta=\frac{1}{4n}$, so that   
 $\underline U(x)$ satisfies  assumption $(G)$.  
We notice that $\Omega_{n}$ is of class $C^2$ for $n\ge N_0$, for certain $N_0$ large. Next we show that $u_{n}$ is a sub-solution of
(\ref{3.2.7}) in $\Omega_{n+1}$. In fact, since $u_{n}$
is the solution of (\ref{3.2.7}) in $\Omega_{n}$ and
$\underline U$ is a sub-solution of (\ref{3.2.7}) in
$\Omega_{n}$, by Theorem \ref{comparison},
$$u_n \ge\underline U\ \ \mbox{in}\ \ \Omega_{n}.$$ 
Additionally, $u_n=\underline U\ \ \mbox{in}\ \
\Omega_{n}^c$. Then, for
$x\in\Omega_{{n+1}}\setminus\Omega_{n}$, we have
\begin{eqnarray*}
(-\Delta)^{\alpha}
u_n(x)=-\frac12\int_{\R^N}\frac{\delta(u_n,x,y)}{|y|^{N+2\alpha}}dy
\leq(-\Delta)^{\alpha}
\underline U(x),
\end{eqnarray*}
so that $u_{n}$ is a sub-solution of (\ref{3.2.7}) in
$\Omega_{n+1}$. From here and since $u_{n+1}$ is the solution of
(\ref{3.2.7}) in $\Omega_{n+1}$ and $\bar U$ is a
super-solution of (\ref{3.2.7}) in $\Omega_{n+1}$,  
by Theorem \ref{comparison}, we have
$
u_{n}\leq u_{n+1}\leq \bar U$ in $\Omega_{n+1}.$
Therefore, for any $n\ge N_0$, $$\underline U\leq u_{n}\leq
u_{n+1}\leq \bar U\ \ \mbox{in}\ \ \Omega.
$$
Then we can  define the function $u$ as
$$
u(x)=\lim_{n\to+\infty}u_n(x),\ x\in \Omega\ \ \mbox{and}\ \
u(x)=0,\ x\in \bar\Omega^c
$$ 
and we have
$$
\underline U(x)\leq u(x)\leq\bar U(x),\  x\in \Omega. 
$$
Since  $\underline U$ and $\bar U$ belong to 
$L^1_\omega(\R^N)$, we see that   $u_n\to u\ \
\mbox{in}\ \  L^1_\omega(\R^N),$ as $n\to\infty$.
Now we repeat the  arguments of the proof of  Theorem \ref{th 3.2.3} to find that u is a classical solution of \equ{eteo1}.
Finally, if $f$ is positive we easily find that  $u$ is positive, again by a contradiction argument. \hfill$\Box$

 \setcounter{equation}{0}
\section{Some estimates}

In order to prove our existence threorems we will use  Theorem \ref{th 3.2.4}, so that it is crucial to have available 
super and sub-solutions to (\ref{3.1.1}).  In this section we provide the basic estimates that will allow to obtain in the next section the necessary super and sub-solutions. 

To this end,  we use appropriate powers of the distance function
$d$ and the main result in this section are the estimates given in Proposition 3.2, that provides the asymptotic behavior of the fractional operator applied to  $d$.

But before going to this estimates, we describe the behavior of the function $C$ defined in \equ{3.1.4}, which is a
$C^2$ defined in  $(-1,2\alpha)$. We have:
\begin{proposition}\label{pr 3.3.2}
 For every  $\alpha\in(0,1)$  there exists a unique $\tau_0(\alpha)\in (-1,0)$ such that
$C(\tau_0(\alpha))=0$ and 
\begin{equation}\label{signotau}
C(\tau)(\tau-\tau_0(\alpha))<0,\quad\mbox{for all}\,\, \tau\in (-1,0)\setminus\{\tau_0(\alpha)\}.
\end{equation}
Moreover, the function $\tau_0$ satisfies
\begin{equation}\label{1.3.2}
\lim_{\alpha\to1^-}\tau_0(\alpha)=0\quad\mbox{and}\quad
\lim_{\alpha\to0^+}\tau_0(\alpha)=-1.
\end{equation}
\end{proposition}

\noindent
{\bf Proof.} We first observe that $C(0)<0$ since the integrand in  \equ{3.1.4} is zero in $(0,1)$ and negative in
$(1,+\infty)$. Next easily see that 
\begin{equation}\label{claim1}
\lim_{\tau\to-1^{+}}C(\tau)=+\infty,
\end{equation}
since, as $\tau$ approaches $-1$, the integrand loses integrability at $0$.
Next we see that  $C(\cdot)$ is strictly convex in $(-1,0)$, since 
$$C'(\tau)=\int_0^{+\infty}\frac{|1-t|^{\tau}\chi_{(0,1)}(t)\log|1-t|+(1+t)^{\tau}\log(1+t)}{t^{1+2\alpha}}dt\
$$and$$
C''(\tau)=\int_0^{+\infty}\frac{|1-t|^{\tau}[\chi_{(0,1)}(t)\log|1-t|]^2+(1+t)^{\tau}[\log(1+t)]^2}{t^{1+2\alpha}}dt>0.$$
The convexity   $C(\cdot)$, $C(0)<0$ and \equ{claim1} allow to conclude the existence and uniqueness of $\tau_0(\alpha)\in (-1,0)$ such that \equ{signotau} holds. To 
 prove the first limit in  (\ref{1.3.2}), we proceed  by contradiction, assuming that
for  $\{\alpha_n\}$ converging to $1$ and $\tau_0\in (-1,0)$ such that
$$\tau_0(\alpha_n)\leq \tau_0<0.$$
Then, for a constant $c_1>0$ we have
$$
\lim_{\alpha_n\to1^{-}}\int^{\frac12}_0\frac{(1-t)^{\tau_0(\alpha_n)}+(1+t)^{\tau_0(\alpha_n)}-2}{t^{1+2\alpha_n}}dt\geq
c_1\lim_{\alpha_n\to1^{-}}\int_{0}^{\frac12}t^{1-2\alpha_n}dt=+\infty
$$
and, for a constant $c_2$ independent of $n$, we have
 \begin{eqnarray*}
\int_{\frac12}^{+\infty}|\frac{\chi_{(0,1)}(t)(1-t)^{\tau_0(\alpha_n)}+(1+t)^{\tau_0(\alpha_n)}-2}{t^{1+2\alpha_n}}|dt&\leq&
c_2,
\end{eqnarray*}
contradicting the fact that 
$C(\tau_0(\alpha_n))=0.$
For the second limit in  (\ref{1.3.2}), we proceed  similarly, assuming that for  $\{\alpha_n\}$ converging to $0$ and 
$\bar \tau_0\in (-1,0)$ such that   
$$\tau_0(\alpha_n)\geq \bar\tau_0>-1.$$
There are positive constants $c_1$ and $c_2$ we have such that  
\begin{eqnarray*}
\int^2_0|\frac{\chi_{0,1}(t)(1-t)^{\tau_0(\alpha_n)}+(1+t)^{\tau_0(\alpha_n)}-2}{t^{1+2\alpha_n}}|dt\le
c_1
\end{eqnarray*}
and
 \begin{eqnarray*}
\lim_{n\to \infty}\int_{2}^{+\infty}\frac{(1+t)^{\tau_0(\alpha_n)}-2}{t^{1+2\alpha_n}}dt&\leq&-c_2\lim_{n\to \infty}\int_2^{+\infty}\frac1{t^{1+2\alpha_n}}dt=
-\infty,
\end{eqnarray*}
contradicting again that $C(\tau_0(\alpha_n))=0.$
\hfill$\Box$\\

Next we prove our main result in this section.  We assume that $\delta>0$ is  such that the distance function
$d(\cdot)$ is of class $C^2$ in $A_\delta=\{x\in \Omega,d(x)<\delta\}$ and we
define
\begin{equation}\label{3.3.1}
V_\tau(x)=\left\{ \arraycolsep=1pt
\begin{array}{lll}
 l(x),\ \ \ \ &
x\in \Omega\setminus A_\delta,\\[2mm]
 d(x)^{\tau},\ & x\in A_\delta,\\[2mm]
0,\ &x\in\Omega^c,
\end{array}
\right.
\end{equation}
 where  $\tau$ is a parameter in $(-1,0)$ and the function $l$  is  positive  such that
$V_\tau$ is $C^2$ in $\Omega$.  We have the following 
\begin{proposition}\label{lm 3.3.1}
Assume $\Omega$ is a bounded, open subset of $\R^N$ with a $C^2$ boundary and let  $\alpha\in(0,1)$. Then 
 there exists $\delta_1\in (0,\delta)$ and a constant $C>1$ such that: 
\\
 $(i)$  If $\tau\in(-1,\tau_0(\alpha))$, then 
$$
\frac1Cd(x)^{\tau-2\alpha }\leq-(-\Delta)^{\alpha}V_\tau(x)\leq
Cd(x)^{\tau-2\alpha },\ \ \mbox{for all}\,\, x\in A_{\delta_1}.$$ 
$(ii)$\ If 
 $\tau\in(\tau_0(\alpha),0)$, then 
$$
\frac1Cd(x)^{\tau-2\alpha }\leq(-\Delta)^{\alpha}V_\tau(x)\leq
Cd(x)^{\tau-2\alpha },\ \ \mbox{for all}\,\, x\in A_{\delta_1}.$$  
$(iii)$\ If 
$\tau=\tau_0(\alpha)$, then 
 $$|(-\Delta)^{\alpha}V_\tau(x)|\leq C
d(x)^{\min\{\tau_0(\alpha),2\tau_0(\alpha)-2\alpha+1\}} , \ \ \mbox{for all}\,\,x\in
A_{\delta_1}.$$
\end{proposition}

\noindent 
{\bf Proof.}
By compactness we prove that the corresponding inequality holds in a neighborhood of any point $\bar x\in\partial\Omega$ and  without loss of generality we may assume that $\bar
x=0$. For  a given $0<\eta\le \delta$, we define
$$
Q_\eta=\{z=(z_1,z')\in\R\times\R^{N-1}, |z_1|<\eta,|z'|<\eta\}$$
and  $Q_\eta^+=\{z\in
Q_\eta,z_1>0\}$.
Let $\varphi:\R^{N-1}\to\R$  be a $C^2$
function such that $(z_1, z')\in \Omega\cap Q_\eta$ if and only if $z_1\in (\varphi(z'),\eta)$ and 
 moreover,   $(\varphi(z'),z')\in\partial\Omega$ for all $|z'|<\eta$. We further assume that  $(-1,0,\cdot\cdot\cdot,0)$ is the outer normal
vector of $\Omega$ at $\bar x$.  

In the proof of our inequalities, we let 
$x=(x_1,0)$, with $x_1\in(0,\eta/4)$, be then a generic point in $A_{\eta/4}$. We observe that $|x-\bar
x|=d(x)=x_1$.
By definition we have
\begin{eqnarray}
-(-\Delta)^\alpha V_\tau(x)=\frac12 \int_{Q_{
\eta}}\frac{\delta(V_\tau,x,y)}{|y|^{N+2\alpha}}dy+\frac12 \int_{\R^N\setminus
Q_{ \eta}}\frac{\delta(V_\tau,x,y)}{|y|^{N+2\alpha}}dy \label{cotadelta}
\end{eqnarray}
and we see that 
\begin{equation}\label{cotadelta1}
|\int_{\R^N\setminus
Q_{ \eta}}\frac{\delta(V_\tau,x,y)}{|y|^{N+2\alpha}}dy|
\le c|x|^\tau,
\end{equation}
where the constant $c$ is independent of $x$. Thus we only need to study the asymptotic behavior of the first integral, that from now on we denote by $E(x_1)$.

Our first goal is to get a lower bound for $E(x_1)$. For that purpose we first notice that, since  $\tau\in(-1,0)$, we have that
\begin{equation}\label{1.3.5}
d(z)^\tau\ge |z_1-\varphi(z')|^\tau,\quad\mbox{for all}\quad  z\in Q_\delta\cap\Omega.
\end{equation}

 Now we assume that  $0<\eta\le \delta/2$,  then for all
$y\in Q_{\eta}$ we have 
 $x\pm y\in Q_\delta $. Thus  $x\pm y\in \Omega\cap Q_\delta$ if and only if
$\varphi(\pm y')<x_1\pm y_1<\delta$ and $|y'|<\delta$.
 Then, by  (\ref{1.3.5}),
we have  that 
\begin{eqnarray}\label{DES1}
\quad V_\tau(x+ y)= d(x+ y)^\tau\ge[x_1+ y_1-\varphi(y')]^\tau,\quad x+y\in Q_{\delta}\cap\Omega 
\end{eqnarray}
and 
\begin{eqnarray}\label{DES2} 
\quad V_\tau(x- y)=d(x- y)^\tau\ge[x_1- y_1-\varphi( -y')]^\tau, \quad x-y\in Q_{\delta}\cap\Omega. 
\end{eqnarray}
  On the other
side, for $y\in Q_{\eta}$ we have that if $x\pm y\in Q_\delta\cap\Omega^c$ then, by definition of $V_\tau$, we have
$V_\tau(x\pm y)=0.$
Now, for  $y\in Q_\eta$ we define the intervals 
\begin{equation}\label{I1}I_+=(\varphi( y')-x_1,\eta-x_1)\quad\mbox{and} \quad I_-=(x_1-\eta,x_1-\varphi(
-y'))
\end{equation}
and the functions
\begin{eqnarray}
I(y)&=&\chi_{I_+}(y_1)|x_1+y_1-\varphi(y')|^{\tau}+\chi_{I_-}(y_1)|x_1-y_1-\varphi(-y')|^{\tau}-2x_1^\tau,~~~~~\nonumber\\
J(y_1)&=&\chi_{(x_1-\eta,x_1)}(y_1)|x_1-y_1|^{\tau}+\chi_{(-x_1,\eta-x_1)}(y_1)|x_1+y_1|^{\tau}-2x_1^\tau, \nonumber \\
I_1(y)&=&\{ \chi_{I_+}(y_1) -\chi_{(-x_1,\eta-x_1)}(y_1)\} |x_1+y_1|^\tau,
 \nonumber\\
 \nonumber I_2(y)&=&
 \chi_{I_+}(y_1) (|x_1+y_1-\varphi(y')|^\tau - 
 |x_1+y_1|^\tau\}), 
 \end{eqnarray}
where $\chi_A$ denotes the characteristic function of the set $A$. Then, using these definitions and inequalities \equ{DES1} and \equ{DES2}, we have that
\begin{eqnarray}\label{Ex1}
\quad E(x_1)\ge \int_{Q_{\eta}} \frac{I(y)}{|y|^{N+2\alpha}}dy=\int_{Q_{\eta}} \frac{J(y_1)}{|y|^{N+2\alpha}}dy+E_1(x_1)+E_2(x_1),
\end{eqnarray}
where 
\begin{equation}\label{Ei}
E_i(x_1)=\int_{Q_{\eta}} \frac{I_i(y)+I_{-i}(y)}{|y|^{N+2\alpha}}dy,\quad i=1,2.
\end{equation}
Here we have considered that 
$$
I_{-1}(y)=\{ \chi_{I_-}(y_1) -\chi_{(x_1-\eta,x_1)}(y_1)\} |x_1-y_1|^\tau
$$
and
$$
 I_{-2}(y)=
 \chi_{I_-}(y_1) (|x_1-y_1-\varphi(-y')|^\tau - 
 |x_1-y_1|^\tau\}),
$$
 for  $y=(y_1,y')\in\R^N$.
We start studying  the first integral in the right hand side in \equ{Ex1}. Changing variables we see that
$$
\int_{Q_{\eta}}\frac{J(y_1)}{|y|^{N+2\alpha}}dy
=x_1^{\tau-2\alpha}\int_{Q_{\frac{\eta}{x_1}}}    \frac{J(x_1z_1)x_1^{-\tau}}{|z|^{N+2\alpha}}dz= 2x_1^{\tau-2\alpha}(R_1-R_2),
$$
where
$$
R_1=\int_{Q_{\frac{\eta}{x_1}}^+}\frac{\chi_{(0,1)}(z_1)|1-z_1|^{\tau}
+(1+z_1)^{\tau}-2}{|z|^{N+2\alpha}}dz
$$
and
$$
R_2=\int_{Q_{\frac{\eta}{x_1}}^+}
\frac{\chi_{(\frac{\eta}{x_1}-1,\frac{\eta}{x_1})}(z_1)(1+z_1)^{\tau}}{|z|^{N+2\alpha}}dz.
$$
Next we estimate these last two integrals. For $R_1$ we see that, for appropriate positive constants $c_1$ and $c_2$ 
 \begin{eqnarray*}
&&\int_{\R^N_+}\frac{\chi_{(0,1)}(z_1)|1-z_1|^{\tau}+(1+z_1)^{\tau}-2}{|z|^{N+2\alpha}}dz\\
&=&\int_0^{+\infty}\frac{\chi_{(0,1)}(z_1)|1-z_1|^{\tau}+(1+z_1)^{\tau}-2}{z_1^{1+2\alpha}}dz_1
\int_{\R^{N-1}}\frac{1}{(|z'|^2+1)^{\frac{N+2\alpha}{2}}}dz'
\\
&=&c_1\,C(\tau)
\end{eqnarray*}
and 
 \begin{eqnarray*}
 \int_{(Q_{\frac{\eta}{x_1}}^+)^c}\frac{\chi_{(0,1)}(z_1)|1-z_1|^{\tau}
+(1+z_1)^{\tau}-2}{|z|^{N+2\alpha}}dz =-c_2\,x_1^{2\alpha}(1+ o(1)).
\end{eqnarray*} 
Consequently we have, for some constant  $c$ that
 \begin{eqnarray}\label{R1}
 R_1=c_1(C(\tau)+cx_1^{2\alpha}+ o(x_1^{2\alpha})).
\end{eqnarray} 
For $R_2$ we have that 
\begin{eqnarray} \quad R_2\label{R2}
&=&\int_{\frac{\eta}{x_1}-1}^{\frac{\eta}{x_1}}\frac{(1+z_1)^{\tau}}{z_1^{1+2\alpha}}
\int_{B_{\frac{\eta}{x_1}}}\frac1{(1+|z'|^2)^{\frac{N+2\alpha}2}}dz'dz_1\le 
c_3x_1^{2\alpha-\tau+1},
\end{eqnarray}
where $c_3>0$. Here and in what follows we denote by $B_\sigma$ the ball of radius $\sigma$ in $\R^{N-1}$.
From \equ{R1} and \equ{R2} we then conclude that
\begin{eqnarray}\label{EQ}
\int_{Q_{\eta}}\frac{J(y_1)}{|y|^{N+2\alpha}}dy=c_1x_1^{\tau-2\alpha}(C(\tau)+cx_1^{2\alpha}+ o(x_1^{2\alpha})).
\end{eqnarray}

Continuing with our analysis we estimate $E_1(x_1)$. 
We only consider the 
 term $I_1(y)$, since the estimate for $I_1(-y)$ is similar. 
We have
\begin{eqnarray*}
 \int_{Q_{\eta}} \frac{I_1(y)}{|y|^{N+2\alpha}}dy
=
-\int_{B_\eta}\int^{\varphi(y')-x_1}_{-x_1}\frac{|x_1+y_1|^\tau}{|y|^{N+2\alpha}}dy_1dy'
 = -x_1^{\tau-2\alpha}F_1(x_1),\label{F1}
\end{eqnarray*}
where
\begin{eqnarray}\label{intF1}
F_1(x_1)= \int_{B_\frac{\eta}{x_1}}\int^{\frac{\varphi(x_1z')}{x_1}}_{0}
\frac{|z_1|^\tau}{((z_1-1)^2+|z'|^2)^{(N+2\alpha)/2}}dz_1dz'.
\end{eqnarray}
In what follows we write $\varphi_-(y')=\min\{\varphi(y'),0\}$ and $\varphi_+(y')=\varphi(y')-\varphi_-(y')$. Next we see that assuming that $0\le \varphi_+(y')\le C|y'|^2$ for $|y'|\le \eta$, for given $(z_1,z')$ satisfying $0\le z_1\le \frac{\varphi_+(x_1 z')}{x_1}$ and $|z'|\le\frac{\eta}{x_1}$ then 
\begin{eqnarray}
\label{cota}(1-z_1)^2+|z'|^2\ge \frac{1}{4}(1+|z'|^2),
\end{eqnarray}
if we assume $\eta$ small enough. Thus
\begin{eqnarray*}
F_1(x_1)&\le&C \int_{B_\frac{\eta}{x_1}}\int^{\frac{\varphi_+(x_1z')}{x_1}}_{0}
\frac{|z_1|^\tau}{(1+|z'|^2)^{(N+2\alpha)/2}}dz_1dz'\\
&\le&  Cx_1^{\tau+1}   \int_{B_\frac{\eta}{x_1}} \frac{|z'|^{2(\tau+1)}}{(1+|z'|^2)^{(N+2\alpha)/2}}dz' \\
&\le &Cx_1^{\tau+1}(x_1^{-2\tau+2\alpha-1}+1)\le Cx_1^{\min\{\tau+1,2\alpha-\tau\}}.
\end{eqnarray*}
Thus we have obtained
\begin{equation}\label{E1(x1)}
E_1(x_1)\ge -Cx_1^{\tau-2\alpha}x_1^{\min\{\tau+1,2\alpha-\tau\}}.
\end{equation}
We continue with 
the estimate  of $E_2(x_1)$.   
As before we 
only consider the term $I_2(y)$,
\begin{eqnarray}
\int_{Q_{\eta}} \frac{I_2(y)}{|y|^{N+2\alpha}}dy&=&\int_{B_\eta}\int_{\varphi(y')-x_1}^{\eta-x_1}\frac{|x_1+y_1-\varphi(y')|^{\tau}-|x_1+y_1|^{\tau}}
{(y_1^2+|y'|^2)^{\frac{N+2\alpha}2}}dy_1dy' \nonumber \\
&\ge&\int_{B_\eta}\int_{\varphi_-(y')-x_1}^{\eta-x_1}\frac{|x_1+y_1-\varphi_-(y')|^{\tau}-|x_1+y_1|^{\tau}}
{(y_1^2+|y'|^2)^{\frac{N+2\alpha}2}}dy_1dy'\nonumber \\
&=&\int_{B_\eta}\int_{\varphi_-(y')}^{\eta}\frac{|z_1-\varphi_-(y')|^{\tau}-|z_1|^{\tau}}
{((z_1-x_1)^2+|y'|^2)^{\frac{N+2\alpha}2}}dz_1dy'\nonumber \\
&\ge&\int_{B_\eta}\int_{0}^{\eta}\frac{|z_1-\varphi_-(y')|^{\tau}-|z_1|^{\tau}}
{((z_1-x_1)^2+|y'|^2)^{\frac{N+2\alpha}2}}dz_1dy' \nonumber \\
& &+ \int_{B_\eta}\int_{\varphi_-(y')}^{0}\frac{-|z_1|^{\tau}}
{((z_1-x_1)^2+|y'|^2)^{\frac{N+2\alpha}2}}dz_1dy'\nonumber\\ &=&
E_{21}(x_1)+E_{22}(x_1).\label{E2cota1}
\end{eqnarray}
We observe that $E_{22}(x_1)$ is similar to $F_1(x_1)$. In order to estimate  $E_{21}(x_1)$ we use integration by parts
\begin{eqnarray*}
E_{21}(x_1)&=&
\frac1{\tau+1}\int_{B_\eta} \left\{\frac{(\eta-\varphi_-(y'))^{\tau+1}-\eta^{\tau+1}}
{((\eta-x_1)^2+|y'|^2)^{\frac{N+2\alpha}2}}    -\frac{(-\varphi_-(y'))^{\tau +1}}{(x_1^2+|y'|^2)^{\frac{N+2\alpha}{2}}} \right\}  dy'\\
&+&
\frac{N+2\alpha}{\tau+1}\int_{B_\eta}\int_{0}^{\eta}\frac{(z_1-\varphi_-(y'))^{\tau+1}-(z_1)^{\tau+1} }
{((z_1-x_1)^2+|y'|^2)^{\frac{N+2\alpha}2+1}}(z_1-x_1)dz_1dy'\\
&=& A_1+A_2.
\end{eqnarray*}
For the first integral we have
\begin{eqnarray*}
A_1&\ge& 
\frac1{\tau+1}\int_{B_\eta} \left\{\frac{-\eta^{\tau+1}}
{((\eta-x_1)^2+|y'|^2)^{\frac{N+2\alpha}2}}    -\frac{(-\varphi_-(y'))^{\tau +1}}{(x_1^2+|y'|^2)^{\frac{N+2\alpha}{2}}} \right\}  dy'
\\
&\ge &
-C(\eta)-C\int_{B_\eta}\frac{|y'|^{2\tau+2}}{(x_1^{2}+|y'|^2)^{\frac{N+2\alpha}2}}dy'
\ge -Cx_1^{\tau-2\alpha+\tau+1} -C.
\end{eqnarray*}
For the second integral, since $\tau\in(-1,0)$ and
$(z_1-\varphi_-(y'))^{\tau+1} -|z_1|^{\tau+1}>0$,  we have that
\begin{eqnarray}\label{q 3.24}
A_2&\ge&\frac{N+2\alpha}{\tau+1}\int_{B_\eta}\int^{x_1}_0\frac{(z_1-\varphi_-(y'))^{\tau+1}
-|z_1|^{\tau+1}}{((z_1-x_1)^2+|y'|^2)^{\frac{N+2\alpha}2+1}}(z_1-x_1)dz_1dy'
\nonumber 
\\[4mm]&\ge&\frac{N+2\alpha}{(\tau+1)^2}\int_{B_\eta}\int^{x_1}_0\frac{-\varphi_-(y')z_1^\tau}{((z_1-x_1)^2+|y'|^2)^{\frac{N+2\alpha}2+1}}(z_1-x_1)dz_1dy'
\nonumber
\\[4mm]&\ge&C_3x_1^{2\tau-2\alpha+1}\int_{B_{\eta/x_1}}\int^{1}_0\frac{|z'|^2z_1^\tau}{((z_1-1)^2+|z'|^2)^{\frac{N+2\alpha}2+1}}(z_1-1)dz_1dz'
\nonumber
\\[4mm]&\ge&-C_4x_1^{2\tau-2\alpha+1},
\end{eqnarray}
where $C_3,C_4>0$ independent of $x_1$ and the second inequality
used $a=z_1$ and $b=-\varphi_-(y')$ in the fact that
$(a+b)^{\tau+1}-a^{\tau+1}\le\frac{a^\tau b}{\tau+1}$ for
$a>0,b\ge0$.

 Thus, we have obtained
\begin{equation}\label{E2(x1)}
E_2(x_1)\ge -Cx_1^{\tau-2\alpha}x_1^{\min\{\tau+1,2\alpha-\tau\}}.
\end{equation}

The next step is to obtain the other inequality for $E(x_1)$. By choosing $\delta$ smaller if necessary, we can prove that
\begin{lema}\label{lemaest}
Under the regularity conditions on the boundary and with the arrangements given at the beginning of the proof, there is $\eta>0$ and $C>0$ such that
$$
d(z)\ge (z_1-\varphi(z'))(1-C|z'|^2)\quad \mbox{for all } (z_1,z')\in \Omega\cap Q_\eta.
$$
\end{lema}

\noindent{\bf Proof.} Since $\varphi$ is $C^2$ and $\nabla \varphi(0)=0$,  there exist $\eta_1\in(0,1/8)$ small and $C_1>0$  such that
$C_1\eta_1<1/4$ and
\begin{equation}\label{e 1}
 |\varphi(y')|<C_1|y'|^2,\quad |\nabla \varphi(y')|\le C_1|y'|,\quad \forall\ y'\in B_{\eta_1}.
\end{equation}
Choosing $\eta_2\in(0,\eta_1)$ such that   for any $z=(z_1,z')\in Q_{\eta_2}\cap\Omega$,
there exists $y'$ satisfying $(\varphi(y'),y')\in\partial\Omega\cap Q_{\eta_1}$ and
$d(z)=|z-(\varphi(y'),y')|$.

We observe that $y'$ mentioned above, is the minimizer of
$$H(z')=(z_1-\varphi(z'))^2+|z'-y'|^2,\quad |z'|<\eta_1,$$
then
$$-(z_1-\varphi(y'))\nabla\varphi(y')+(z'-y')=0,$$
which, together with (\ref{e 1}) implies that
\begin{eqnarray*}
|y'|-|z'|&\le&|z'-y'|=|(z_1-\varphi(y'))\nabla\varphi(y')|
\le(|z_1|+C_1|y'|^2)|\nabla\varphi(y')|
\\ &\le& C_1(\eta_2+C_1\eta_1^2) |y'|
\le 2C_1\eta_1 |y'|<\frac12|y'|.
\end{eqnarray*}
Then
\begin{equation}\label{e 2}
|y'|\le 2|z'|.
\end{equation}

Denote the points $z,(\varphi(y'),y'),(\varphi(z'),z') $ by $A,B,C$, respectively, and
let $\theta$ be the angle between the segment $BC$ and the hyper plane with normal vector $e_1=(1,0,...,0)$ and containing $C$.
Then the angle $\angle C=\frac\pi2-\theta.$ Denotes the arc from $B$ to $C$ in the plane $ABC$  by arc$(BC)$. By the geometry, there exists some point
$x=(\varphi(x'),x')\in \mbox{arc}(BC)$ such that line $BC$ parallels the tangent line of arc$(BC)$ at point $x$. Then, from (\ref{e 2}) we have $|x'|\le \max\{|z'|,|y'|\}\le 2|z'|$ and so, from (\ref{e 1}) we obtain
\begin{eqnarray*}
\tan( \theta)= |\frac{y'-z'}{|y'-z'|}\cdot\nabla\varphi(x')|\le|\nabla\varphi(x')|
\le C_1|x'|\le 2C_1|z'|,
\end{eqnarray*}
which implies that for some $C>0$,
\begin{equation}\label{e 3}
\cos(\theta)\ge 1-C|z'|^2.
\end{equation}
Then we complete the proof  using Sine Theorem and (\ref{e 3})
\begin{eqnarray*}
d(z)&=&\frac{\sin(\angle C)}{\sin(\angle B)}(z_1-\varphi(z'))
\ge (z_1-\varphi(z'))\sin(\frac\pi2-\theta)
\\&=&(z_1-\varphi(z'))\cos(\theta)
\ge (z_1-\varphi(z'))(1-C|z'|^2). \qquad \Box
\end{eqnarray*}

\medskip

From this lemma,  by making $C$ and $\eta$ smaller if necessary we obtain that
\begin{equation}\label{eq 14}
d^\tau(z)\le (z_1-\varphi(z'))^\tau (1+C|z'|^2)\quad \mbox{for all }  z\in \Omega\cap Q_\eta.
\end{equation}
With  $x=(x_1,0)$ satisfying $x_1\in(0,\eta/4)$ as at the beginning of the proof, we have that 
 $d(x)=x_1$ and for any $y\in Q_{\eta}$ we see that
 $x\pm y\in Q_\delta.$
We also see that $x\pm y\in \Omega\cap Q_{\delta} $ if and only if
$\varphi(\pm y')<x_1\pm y_1<\delta$ and $|y'|<\delta$. Then, for $x\pm y\in
\Omega\cap Q_{\eta} $, by (\ref{eq 14}) we have,
\begin{eqnarray}\label{B1}
V_\tau(x\pm y)=d(x\pm y)^\tau\le (x_1\pm
y_1-\varphi(\pm y'))^\tau(1+C|y'|^2).
\end{eqnarray}
%
%
%
%
For $y\in Q_\eta$, we define
$$
I_3(y)=C|y'|^2 \chi_{I_+}(y_1)|x_1+y_1-\varphi(y')|^\tau
$$
and
$$
I_3(-y)=C|y'|^2 \chi_{I_-}(y_1)|x_1-y_1-\varphi(-y')|^\tau,
$$
where $I_+$ and $I_-$ were defined in \equ{I1}. 
Using \equ{B1} as in \equ{Ex1} we find 
\begin{eqnarray}\label{B2}
E(x_1)&=&\int_{Q_\eta}\frac{\delta(V_\tau,x,y)}{|y|^{N+2\alpha}}dy\le \int_{Q_\eta}\frac{I(y)}{|y|^{N+2\alpha}}dy+E_3(x_1)
\nonumber\\
&=&\int_{Q_\eta}\frac{J(y)}{|y|^{N+2\alpha}}dy+E_1(x_1)+E_2(x_1)+E_3(x_1),
\end{eqnarray}
where $E_1$ and $E_2$ were defined in \equ{Ei} and 
\begin{equation}\label{E33}
E_3(x_1)=\int_{Q_\eta}\frac{I_3(y)+I_3(-y)}{|y|^{N+2\alpha}}dy.
\end{equation}
We estimate $E_3(x_1)$ and for that we observe that it is enough to estimate the integral with one of the terms in \equ{E33} (the other is similar), say
\begin{eqnarray}
& &\int_{Q_\eta}\frac{I_3(y)}{|y|^{N+2\alpha}}dy=\int_{B_\eta}\int_{\varphi(y')-x_1}^{\eta-x_1} \frac{C|y'|^2 |x_1+y_1-\varphi(y')|^\tau}{|y|^{N+2\alpha}}dy_1dy'\nonumber\\
\quad\quad  &=& Cx_1^{\tau-2\alpha+2}\int_{B_{\frac{\eta}{x_1}}}\int_{\frac{\varphi(x_1z')}{x_1}}^{\frac{\eta}{x_1}} \frac{|z'|^2 |z_1-\frac{\varphi(x_1z')}{x_1}|^\tau}{((z_1-1)^2+|z'|^2)^{(N+2\alpha)/2}}dz_1dz'\nonumber\\
\quad \quad &=&Cx_1^{\tau-2\alpha+2}(A_1+A_2), \label{A10}
\end{eqnarray}
where $A_1$ and $A_2$ are integrals over properly chosen subdomains, estimated separately.
\begin{eqnarray}
A_1&=&\int_{B_{\frac{\eta}{x_1}}}\int_{\frac{\varphi(x_1z')}{x_1}}^{\frac{\varphi(x_1z')}{x_1}+\frac{1}{2}} \frac{|z'|^2 |z_1-\frac{\varphi(x_1z')}{x_1}|^\tau}{((z_1-1)^2+|z'|^2)^{(N+2\alpha)/2}}dz_1dz' \nonumber\\
&\le& 
\frac{c}{(\tau+1)2^{\tau+1}}\int_{B_{\frac{\eta}{x_1}}}\frac{|z'|^2 }{(1+|z'|^2)^{(N+2\alpha)/2}}dz'  \label{A11}   \\
&\le&  c'\left(\frac{\eta}{x_1}\right)^{-2\alpha+1}.\label{A12}
\end{eqnarray}
The inequality in \equ{A11} is obtained noticing that the ball $B((1,0),1/2)$ in $R^N$  does not touch the band
$$
\{(z_1,z')\,/\, |z'|\le \eta,  \frac{\varphi(x_1z')}{x_1}\le z_1\le  \frac{\varphi(x_1z')}{x_1}+1/2\}
$$
if $x_1$ is small enough, and so $(z_1-1)^2+|z'|^2
 \ge \frac18+ \frac12|z'|^2$. Then simple integration gives the next term. 
 Next we estimate 
 $A_2$
 \begin{eqnarray}
A_2&=&\int_{B_{\frac{\eta}{x_1}}}\int^{\frac{\eta}{x_1}}_{\frac{\varphi(x_1z')}{x_1}+\frac{1}{2}} 
 \frac{|z'|^2 |z_1-\frac{\varphi(x_1z')}{x_1}|^\tau}{((z_1-1)^2+|z'|^2)^{(N+2\alpha)/2}}dz_1dz' \nonumber\\
 &\le& \frac1{2^\tau} \int_{B_{\frac{\eta}{x_1}}}\int^{\frac{\eta}{x_1}}_{\frac{\varphi(x_1z')}{x_1}+\frac{1}{2}} 
 \frac{|z'|^2 }{((z_1-1)^2+|z'|^2)^{(N+2\alpha)/2}}dz_1dz' \nonumber\\
 &\le&  c'\left(\frac{\eta}{x_1}\right)^{-2\alpha+2}.\label{A13} 
 \end{eqnarray}
Putting together \equ{A10},   \equ{A12}, \equ{A13}  and \equ{E33} we obtain
\begin{equation}
E_3(x_1)=\int_{Q_\eta}\frac{(I_3(y)+I_3(-y))}{|y|^{N+2\alpha}}dy\le  cx_1^{\tau}.\label{I333}
\end{equation}
From \equ{F1},  but using the other inequality for $F_1$, that is,
$$
F_1(x_1)\ge C \int_{B_\frac{\eta}{x_1}}\int^{\frac{\varphi_-(x_1z')}{x_1}}_{0}
\frac{|z_1|^\tau}{(1+|z'|^2)^{(N+2\alpha)/2}}dz_1dz'
$$
 and arguing similarly we obtain as in \equ{E1(x1)}
\begin{equation}\label{E1(x1)+}
E_1(x_1)\le Cx_1^{\tau-2\alpha}x_1^{\min\{\tau+1,2\alpha\}}.
\end{equation}
Then we look at $E_2(x_1)$ and, as  in \equ{E2cota1}, we  
only consider the term $I_2(y)$: 
\begin{eqnarray*}
\int_{Q_{\eta}} \frac{I_2(y)}{|y|^{N+2\alpha}}dy&
\le&\int_{B_\eta}\int_{\varphi_+(y')}^{\eta}\frac{|z_1-\varphi_+(y')|^{\tau}-|z_1|^{\tau}}
{((z_1-x_1)^2+|y'|^2)^{\frac{N+2\alpha}2}}dz_1dy' =\tilde E_{21}(x_1).
\end{eqnarray*}
In order to estimate  $\tilde E_{21}(x_1)$ we use integration by parts
\begin{eqnarray*}
& &\tilde E_{21}(x_1)=\\
& &\! \! \! \! \frac1{\tau+1}\int_{B_\eta} \left\{\frac{(\eta-\varphi_+(y'))^{\tau+1}-\eta^{\tau+1}}
{((\eta-x_1)^2+|y'|^2)^{\frac{N+2\alpha}2}}    -\frac{(\varphi_+(y'))^{\tau +1}}{((\varphi_+(y')-x_1)^2+|y'|^2)^{\frac{N+2\alpha}{2}}} \right\}  dy'\\
& & +
\frac{N+2\alpha}{\tau+1}\int_{B_\eta}\int_{\varphi_+(y')}^{\eta}\frac{(z_1-\varphi_+(y'))^{\tau+1}-z_1^{\tau+1} }
{((z_1-x_1)^2+|y'|^2)^{\frac{N+2\alpha}2+1}}(z_1-x_1)dz_1dy'\\
& &\le 
\frac{N+2\alpha}{\tau+1}\int_{B_\eta}\int_{\min\{\varphi_+(y'),x_1\}}^{x_1}\frac{ (z_1-\varphi_+(y'))^{\tau+1}-z_1^{\tau+1}   }
{((z_1-x_1)^2+|y'|^2)^{\frac{N+2\alpha}2+1}}(z_1-x_1)dz_1dy'.
\end{eqnarray*}
This integral can be estimated in a similar way as $E_{21}$, see  \equ{q 3.24} and the estimates given before.  We then obtain 
\begin{equation}\label{E2(x1)+}
E_2(x_1)\le Cx_1^{2\tau-2\alpha+1}.
\end{equation}
Then we conclude from \equ{cotadelta},  \equ{Ex1},  \equ{EQ},  \equ{E1(x1)},  \equ{E2(x1)},  \equ{B2}, \equ{I333},  \equ{E1(x1)+} and  \equ{E2(x1)+}
   that
\begin{eqnarray}
-(-\Delta)^\alpha V_\tau(x)= Cx_1^{\tau-2\alpha}(C(\tau)+O(x_1^{\min\{\tau +1, 2\alpha\}})),
\end{eqnarray}
where there exists a constant $c>0$ so that
$$
|O(x_1^{\min\{\tau +1, 2\alpha\}})|\le c x_1^{\min\{\tau +1, 2\alpha\}},\qquad \mbox{for all small } x_1>0.
$$
From here, depending on the value of $\tau\in (-1,0)$, conditions (i), (ii) and (iii) follows and the proof of the proposition is complete.$\hfill\Box$

\medskip

We end this section with an estimate we need when dealing with equation \equ{3.1.1} when the external value $g$ is not zero. We have the following proposition 
\begin{proposition}\label{lm 3.3.21}
Assume that $\Omega$ is a bounded, open and $C^2$ domain in $\R^N$.
Assume that  $g\in L^1_\omega(\Omega^c)$. Assume further that there are numbers $\beta\in (-1,0)$, $\eta>0$ and $c>1$ such that
$$\frac1c\leq g(x) d(x)^{-\beta}\leq c,\ \ x\in\bar\Omega^c\ \mbox{and}\ d(x)\leq \eta.$$
 Then there exist
$\eta_1>0$ and $C>1$
 such that  $G$, defined in \equ{G}, satisfies
 \begin{equation}\label{3.3.2}
 \frac1Cd(x)^{\beta-2\alpha }\leq G(x)\leq
Cd(x)^{\beta-2\alpha },\ \ x\in A_{\eta_1}.
\end{equation}
\end{proposition}
\noindent
{\bf Proof.} The proof of this  proposition requires estimates  similar to those in the proof of Proposition \ref{lm 3.3.1} so we omit it. However, the function $C$ used there and defined in \equ{3.1.4},  needs to be replaced here by $\tilde C: (-1,0)\to \R$ given by
$$
\tilde C(\beta)=\int_{1}^\infty\frac{|t-1|^\beta}{t^{1+2\alpha}}dt.
$$
We observe that this function is always positive.   \hfill$\Box$

\section{ Proof of existence results}

In this section, we will give the proof of existence of large
solution to \equ{eteo1}. 
By Theorem  \ref{th 3.2.4} we only need to find ordered super and sub-solution, denoted by $U$ and $W$, for \equ{eteo1} under
our various assumptions.
We begin with a simple lemma that reduce the problem to find them only in $A_\delta$.

\begin{lemma}
Let $U$ and $W$  be classical ordered super and sub-solution  of \equ{eteo1}  in the sub-domain  $A_\delta$.
Then there exists $\lambda$ large such that
$U_\lambda=U-\lambda\bar V$ and  $ W_\lambda=W+\lambda\bar V$,
where $\bar V$ is the solution of  \equ{3.2.V}, with ${\mathcal{O}}=\Omega$,  are ordered super and sub-solution  of \equ{eteo1}.
\end{lemma}

\noindent {\bf Proof.} 
Notice that by negativity $\bar V$ in $\Omega$, we have that $U_\lambda\geq U$ and $W_\lambda\leq W$, so they are still ordered in $A_\delta$. In addition $U_\lambda$ satisfies $$(-\Delta)^\alpha U_\lambda+|U_\lambda|^{p-1}U_\lambda-f(x)\ge (-\Delta)^\alpha U+|U|^{p-1}U-f(x)+\lambda>0,\quad\mbox{in}\quad\Omega.$$
This inequality holds because of our  assumption in $A_\delta$, the fact that $ (-\Delta)^\alpha U+|U|^{p-1}U-f(x)$ is continuous in $\Omega\setminus{A_\delta}$ and
by  taking $\lambda$ large enough. 

By the same type of arguments we find the $W_\lambda$ is a sub-solution  of the first equation in \equ{eteo1} and we complete the proof.
$\Box$

\medskip

Now we are in position to prove our existence results  that we already reduced to find  ordered super and sub-solution  of \equ{eteo1}  with the first equation in
 $A_\delta$ with the desired asymptotic behavior.
 
\noindent {\bf Proof of Theorem \ref{th 3.1.1} (Existence).} Define
\begin{equation}\label{UUWW}
U_\mu(x)=\mu V_\tau(x) \quad \mbox{and}\quad \ W_{\mu}(x)=\mu V_\tau(x),
\end{equation}
with $\tau=-\frac{2\alpha}{p-1}$. We observe  that $\tau=-\frac{2\alpha}{p-1}\in (-1,\tau_0(\alpha))$
 and $\tau p=\tau-2\alpha $, 
Then by Proposition  \ref{lm 3.3.1} and $(H2)$ we find that for $x \in A_\delta$ and $\delta>0$ small
\begin{eqnarray*}
(-\Delta)^\alpha U_\mu(x)+U^p_\mu(x)-f(x)\geq -C\mu
d(x)^{\tau-2\alpha}+\mu^p d(x)^{\tau
p}-Cd(x)^{\tau p},
\end{eqnarray*} for some $C>0$. Then there exists a large $\mu>0$ such that $U_\mu$ is a 
super-solution  of \equ{eteo1}  with the first equation in  $A_\delta$ with the desired asymptotic behavior.
Now  by Proposition  \ref{lm 3.3.1} we have that for $x \in A_\delta$ and $\delta>0$ small
\begin{eqnarray*}
(-\Delta)^\alpha W_\mu(x)+W^p_\mu(x)-f(x)\leq -\frac{\mu}{C}
d(x)^{\tau-2\alpha}+\mu^p d(x)^{\tau p}-f(x)\leq 0,
\end{eqnarray*}
in the last inequality we have used $(H2)$ and $\mu>0$ small.  Then, by Theorem \ref{th 3.2.4} there exists a solution, with the desired asymptotic behavior. $\Box$

\noindent {\bf Proof of Theorem \ref{th 3.1.1} (Special case $\tau=\tau_0(\alpha)$)}.  We define for $t>0$,
\begin{equation}\label{obs3.4.1}
U_{\mu}(x)=t V_{\tau_0(\alpha)}(x)-\mu
V_{\tau_1}(x) \quad \mbox{ and}\quad  W_{\mu}(x)=t
V_{\tau_0(\alpha)}(x)-\mu V_{\tau_1}(x),
\end{equation}
where $\tau_1=\min\{\tau_0(\alpha)p+2\alpha,0\}$.
If $\tau_1=0$, we write $V_0=\chi_\Omega$ and we have
\begin{eqnarray*}
(-\Delta)^\alpha
V_0(x)=\int_{\R^N\setminus\Omega}\frac1{|z-x|^{N+2\alpha}}dz,\quad x\in\Omega.
\end{eqnarray*}
By direct computation, there exists $C>1$ such that
\begin{equation}\label{V_0}
\frac1C d(x)^{-2\alpha}\le(-\Delta)^\alpha V_0(x)\le C
d(x)^{-2\alpha},\quad x\in\Omega.
\end{equation}
We see that $\tau_1\in(\tau_0(\alpha),0]$  and, if $\tau_1<0$,  we have $\tau_1-2\alpha=\tau_0(\alpha)p$
and
$$\tau_1-2\alpha<\min\{\tau_0(\alpha),\tau_0(\alpha)-2\alpha+\tau_0(\alpha)+1\}.$$, u
Then,
by Proposition \ref{lm 3.3.1} and (\ref{V_0}), for $x\in A_{\delta}$, it
follows that
\begin{eqnarray*}
(-\Delta)^\alpha U_{\mu}(x)+|U_{\mu}(x)|^{p-1}U_{\mu}(x)&\geq& -Ct
d(x)^{\min\{\tau_0(\alpha),\tau_0(\alpha)-2\alpha+\tau_0(\alpha)+1\}
}\\&&-C \mu  d(x)^{\tau_1-2\alpha }+t^p d(x)^{\tau_0(\alpha)
p}.\end{eqnarray*} 
Thus, letting $\mu=t^p/(2C)$ if $\tau_1<0$ and $\mu=0$  if $\tau_1=0$,  for a possible smaller $\delta>0$, we obtain
$$(-\Delta)^\alpha U_{\mu}(x)+|U_{\mu}(x)|^{p-1}U_{\mu}(x)\geq 0,\quad x\in
A_{\delta}.$$
For the sub-solution, by Proposition \ref{lm 3.3.1} and (\ref{V_0}), for $x\in A_{\delta}$, we
have
\begin{eqnarray*}(-\Delta)^\alpha
 W_{\mu}(x)+ |W_{\mu}|^{p-1}W_{\mu}(x)&\leq&
Ct
d(x)^{\min\{\tau_0(\alpha),\tau_0(\alpha)-2\alpha+\tau_0(\alpha)+1\}
}\\&&-\frac\mu C d(x)^{\tau_1-2\alpha }+t^p d(x)^{\tau_0(\alpha)
p},\end{eqnarray*} where $C>1$. 
Then, for $\mu\ge2Ct^p$ and a possibly smaller  $\delta>0$
$$(-\Delta)^\alpha  W_{\mu}(x)+|W_{\mu}|^{p-1}W_{\mu}(x)\leq 0,\ x\in
A_{\delta},$$
completing the proof.\hfill $\Box$

\noindent {\bf Proof  of Theorem  \ref{th 3.1.2}.}
We define $U_\mu$ and $W_\mu$ as in \equ{UUWW}.
In the case of a weak source,  we take $\tau=\gamma+2\alpha$ and we observe that $\gamma+2\alpha\ge-\frac{2\alpha }{p-1}>\tau_0(\alpha)$ and $p(\gamma+2\alpha)\geq \gamma$. 
Using Proposition  \ref{lm 3.3.1} and $(H3)$ we find that $U_\mu$  is a super-solution  for $\mu>0$ large (resp. $W_\mu$ 
is a sub-solution for $\mu>0$ small)  of \equ{eteo1}  with the first equation in
 $A_\delta$ for $\delta>0$ small.
In the case of a strong source, we take $\tau=\frac{\gamma}{p}$ and observe that $\gamma <\frac{\gamma}{p}-2\alpha$.
Using Proposition  \ref{lm 3.3.1} we find
$$|(-\Delta)^\alpha U_mu|,|(-\Delta)^\alpha U_mu|\leq C d(x)^{\frac{\gamma}{p}-2\alpha}.$$
By $(H3)$  we find that $U_\mu$  is a super-solution  for $\mu$ large (resp. $W_\mu$ 
is a sub-solution for $\mu$ small)  of \equ{eteo1}  with the first equation in
 $A_\delta$ for $\delta$ small.\hfill $\Box$

\begin{remark}\label{RQQQ}
In order to obtain the above existence results  for classical solution to (\ref{3.1.1}), that is when $g$ is not necessarily zero, we only need use them with $F$ as a right hand side as given in \equ{F}. Here we only need to assume that $g$ satisfies  $(H4)$. 
In fact, as above we find super and sub-solutions  for \equ{eteo1},  with $f$ replaced by $F$.  Then, as in the proof of  Theorem \ref{th 3.2.4},  we find a viscosity solution of \equ{eteo1} and then  $v=u+\tilde g$ is a viscosity solution of (\ref{3.1.1}).  Next we use Theorem 2.6 in  \cite{CS2} and then we use  Theorem \ref{regularity} to obtain that $v$ is a classical solution of (\ref{3.1.1}).\\
\end{remark}

\begin{remark}\label{RQQQ1}\noindent Now we compare Theorem \ref{th 3.1.1} with the result in \cite{FQ0}.
 Let us assume that $f$ and $g$ satisfies hypothesis (F0)-(F2) and (G0)-(G3), respectively, given in \cite{FQ0}. We first 
observe that the function $F$, as defined above, satisfies $(H1)$ thanks to (G0), (G3) and (F0). Next we see
 that $F$ satisfies $(H2)$, since (G2), (F1) and (F2) holds. Here we have to use Proposition \ref{lm 3.3.21}.
  In the range of $p$ given by \equ{pp1}, we then may apply Theorem \ref{th 3.1.1} to obtain existence of a 
  blow-up solution as given in Theorem 1.1 in \cite{FQ0}. We see that the existence is proved here, without 
  assuming hypothesis $(G1)$, thus we generalized this earlier result. Moreover, here we obtain a uniqueness and
   non existence of blow-up solution, if we  further assume hypotheses on $f$ and $g$, guaranteeing hypothesis 
   $(H2^*)$ in Theorem  \ref{th 3.1.1}. The complementary range of $p$ is obtained using Theorem 1.2 for the existence of solutions as given in  Theorem 1.1 in \cite{FQ0} and uniqueness and non-existence as in Theorem 1.3 and 1.4 are truly new results. The hypotheses needed on $g$ to obtain $(H3)$ for the function $F$  are a bit stronger, since we are requiring in $(H3)$ that the explosion rate is the same from above and from below, while in (G2) and (G4) they may be different.  

  \end{remark}

\setcounter{equation}{0}
\section{Proof of uniqueness results}

In this section we prove our uniqueness results, which are given in Theorem \ref{th 3.1.1} and Theorem \ref{th 1.1}. These results are for positive solutions, so we assume that the external source $f$ is non-negative. We  assume that there are two positive solutions $u$ and $v$ of  \equ{eteo1} and then define the set
\begin{equation}\label{4.1.4}
\mathcal{A}=\{x\in  \Omega,\ u(x)>v(x)\}.
\end{equation}
This set is open, $\mathcal{A}\subset \Omega$ and we only need  to prove that $\mathcal{A}=\O,$ to obtain that $u=v$, by interchanging the roles of $u$ and $v$.

We will distinguish three cases, depending on the conditions satisfying $u$ and $v$: 
Case a)  $u$ and $v$  satisfy \equ{pp1} and \equ{3.1.5} (uniqueness part of Theorem \ref{th 3.1.1}),
Case b) $u$ and $v$  $\equ{gamma1}$ and \equ{3.1.10} (weak source in Theorem \ref{th 1.1}) and 
Case c) $u$ and $v$ with $\equ{gamma2}$ -\equ{3.1.11} (strong source in Theorem \ref{th 1.1}).

We start our proof considering an auxiliary function  
\begin{equation}\label{epv}
V(x)=\left\{ \arraycolsep=1pt
\begin{array}{lll}
 c(1-|x|^2)^3,\ \ \ \ &
x\in B_1(0),\\[2mm]
0,\ &x\in B_1^c(0),
\end{array}
\right.
\end{equation}
where the constant $c$ may be chosen so that $V$  satisfies 
\begin{equation}\label{2.3} (-\Delta)^\alpha V(x)\leq 1\quad\mbox{and}\quad  0<V(0)=\max_{x\in\R^N}V(x).
\end{equation}

In order to prove the uniqueness result in the three cases, we need first some preliminary lemmas.
\begin{lemma}\label{lm 3.1}
If
$\mathcal{A}_{k}=\{x\in  \Omega, u(x)-k v(x)> 0\}\not=\O,$
for $k>1$. Then,
\begin{equation}\label{4.1.5}
\partial \mathcal{A}_{k}\cap\partial \Omega\not=\O.
\end{equation}
\end{lemma}
\noindent
{\bf Proof.}  If (\ref{4.1.5}) is not true, there exists 
$\bar x\in \Omega$ such that
$$u(\bar x)-kv(\bar x)=\max_{x\in\R^N}(u-kv)(x)>0,$$
 Then, we have
$$
(-\Delta)^\alpha (u-k v)(\bar x)\geq0,
$$
which contradicts 
\begin{eqnarray*}
(-\Delta)^\alpha (u-k v)(\bar x)&=&-u^p(\bar x)+kv^p(\bar
x)-(k-1)f(\bar x)\\&\leq&-(k^p-k)v^p(\bar x)
<0.\hfill\Box
\end{eqnarray*}

\begin{lemma}\label{lm 3.2}
If  
$\mathcal{A}_{k}\not=\O$,
for $k>1$, then
 \begin{equation}\label{4.1.6}
 \sup_{x\in\Omega}(u-kv)(x)=+\infty.
 \end{equation}
\end{lemma}
\noindent
{\bf Proof.} Assume that  $\bar M=\sup_{x\in \Omega}(u-kv)(x)<+\infty.$ We see that $\bar M>0$ and
there is no point $\bar x \in \Omega$ achieving the supreme of $u-kv$, by the same argument given above. 
Let us consider $x_0\in\mathcal{A}_{k}$,
$r=d(x_0)/2$ and define
 \begin{equation}\label{wk}
 w_k=u-kv\ \ \mbox{in}\ \ \R^N.
 \end{equation}
 Under the conditions of Case a)  and b) (resp. Case c)),
for all $x\in B_r(x_0)\cap\mathcal{A}_{k}$ we have
\begin{equation}\label{4.2.1}
 (-\Delta)^\alpha
 w_k(x)=-u^p(x)+kv^p(x)+(1-k)f(x)\le-K_1r^{\tau-2\alpha}, 
\end{equation}
(resp. $  \leq-K_1r^{\gamma}$). Here we have used that $\tau=-2\alpha/(p-1)$ and, in Case a) \equ{3.1.5} for $v$, in Case b) $(H3)$ and \equ{gamma1} and in Case c)  $(H3)$.
Moreover, in Case a) we have considered $K_1=C(k^p-k)$ and in Cases b) and  c) $K_1=C(k-1)$ for some constant $C$.
Now we define 
$$
w(x)=\frac{2\bar M}{V(0)} V\left (\frac{x-x_0}{r}\right )
$$ 
for $x\in \R^N,$
where $V$ is given in \equ{epv}, and  we see that
\begin{equation}\label{4.1.8}
w(x_0)=2\bar M
\end{equation}
and
\begin{equation}\label{4.1.9}
(-\Delta)^\alpha w\leq \frac{2\bar M}{V(0)} r^{-2\alpha},\ \ \ \
\mbox{in}\ \ B_r(x_0).
\end{equation}
Since $\tau<0$ ($\gamma<-2\alpha$ in the Case c)), by Lemma \ref{lm 3.1}  we can take $x_0\in\mathcal{A}_{k}$ close to
$\partial\Omega$, so that
$$
\frac{2\bar M}{V(0)}\leq K_1r^\tau\
\ \ ( \ \frac{2\bar
M}{V(0)}\leq K_1r^{\gamma+2\alpha},\quad  \mbox{ in Case c))}.
$$ 
From here, combining (\ref{4.2.1}) with
(\ref{4.1.9}), we have that
\begin{eqnarray*}
 (-\Delta)^\alpha (w_k+w)(x)\le0,\ \ \ x\in
 B_r(x_0)\cap\mathcal{A}_{k}.
\end{eqnarray*}
Then, by the Maximum Principle, we obtain
\begin{equation}\label{ewk}
w_k(x_0)+w(x_0)\leq \max\{\bar M, \sup_{x\in B_r(x_0)\cap \mathcal{A}_{k}^c}(w_k+w)\}.
\end{equation}
In case we have
\begin{equation}\label{ewk1}
\bar M<\sup_{x\in B_r(x_0)\cap \mathcal{A}_{k}^c}(w_k+w),
\end{equation}
then
\begin{eqnarray}\label{mpw}
 w(x_0)<(w_k+w)(x_0)&\le &\sup_{x\in B_r(x_0)\cap \mathcal{A}_{k}^c}(w_k+w)(x)\nonumber \\
& \le&   \sup_{x\in
B_r(x_0)\cap \mathcal{A}_{k}^c}w(x) =w(x_0),  
\end{eqnarray}
which is impossible. 
So that \equ{ewk1} is false and then, from  \equ{ewk} we get
$$w(x_0)<w_k(x_0)+w(x_0)\leq \bar M,$$
which is impossible in view of (\ref{4.1.8}), completing  the proof. \hfill
$\Box$

\begin{lemma}\label{lm 3.3}
 There exists a sequence  $\{C_n\}$, with $C_n>0$, satisfying
\begin{equation}\label{4.1.10}
 \lim_{n\to+\infty} C_n=0
\end{equation}
and such that for all $x_0\in
\mathcal{A}_{k}$ and $k>1$ we have
\begin{eqnarray*}
0<\int_{Q_n}\frac{w_k(z)-M_n}{|z-x|^{N+2\alpha}}dz\leq
C_nr^{\tau-2\alpha},\ \ \forall x\in B_r(x_0),
\end{eqnarray*}
where we consider $r=d(x_0)/2$, $Q_n=\{z\in  A_{r/n}\, / \, w_k(z)>M_n\}$ and
$M_n=\max_{x\in\Omega\setminus A_{r/n}}w_k(x).$
\end{lemma}
\noindent
{\bf Proof.}  In Case a): we see that $Q_n\subset   A_{r/n}$ and
$\lim_{n\to+\infty}|Q_n|=0$, so that using \equ{1.7} we directly  obtain
\begin{eqnarray*}
\int_{Q_n}\frac{w_k(z)-M_n}{|z-x|^{N+2\alpha}}dz&\leq&C_0r^{-N-2\alpha}\int_{
 A_{r/n}}d(z)^{\tau}dz \\
&\le&Cr^{-N-2\alpha}\int_0^{r/n}t^{\tau}t^{N-1}dt
\le\frac C{n^{N+\tau}}r^{\tau-2\alpha},
\end{eqnarray*}
where $C$ depends on $C_0$ and $\partial\Omega$. We complete the proof defining $C_n=\frac
C{n^{N+\tau}}$.

In Case b) we argue similarly using \equ{3.1.10} and define $C_n$ as before, while in Case c) we argue similarly using 
  \equ{3.1.11}, but defining
$C_n=\frac {C}{n^{N+\gamma /p}}$. \hfill$\Box$

\medskip

Now we are in a position to prove our non-existence results.

\noindent
{\bf Proof of uniqueness results in Cases a), b) and c).}
We assume that ${\mathcal A}\not=\O$, then there  exists $k>1$ such that ${\mathcal A}_k\not=\O$. 
By Lemma \ref{lm 3.2} 
there exists $x_0 \in {\mathcal A}_k$ such that
$$
w_k(x_0)=\max\{w_k(x)\, /\, x\in \Omega\setminus A_{d(x_0)} \}.
$$
Proceeding as in Lemma \ref{lm 3.2} with the function 
$$
w(x)=\frac{K_1}{2}r^\tau V(\frac{x-x_0}{r})
$$
$$
\mbox{and}\quad  
w(x)=\frac{K_1}{2}r^ {\gamma+2\alpha}V(\frac{x-x_0}{r}), \mbox{ in Case c)}, 
$$
we see that
\begin{equation}\label{ref*}
 (-\Delta)^\alpha (w_k+w)(x)\le -\frac{K_1}{2}r^{\tau-2\alpha} ,\ \ \ x\in
 B_r(x_0)\cap\mathcal{A}_{k}.
\end{equation}
\begin{equation}
\mbox{and}\quad
 (-\Delta)^\alpha (w_k+w)(x) \le -\frac{K_1}{2}r^{\gamma},	\mbox{ in Case c).}
 \end{equation}
With $M_n$, as given in Lemma \ref{lm 3.3},
we define
\begin{equation} \label{gg3.2.5}
\bar w_n(x)=\left\{ \arraycolsep=1pt
\begin{array}{ll}
(w_{k}+w)(x),\ \ \ \ &\mbox{if}\ \
 w_{k}(x)\leq M_n   ,\\[2mm]
M_n,&\mbox{if}\ \ w_{k}(x)> M_n,
\end{array}
\right.
\end{equation}
for $n>1$.
By Lemma \ref{lm 3.3}  we find  $n_0$ such that  
\begin{eqnarray*}
(-\Delta)^\alpha\bar w_{n_0}(x)&=&(-\Delta)^\alpha
(w_{k}+w)(x)+2\int_{
Q_{n_0}}\frac{w_{k}(z)-M_{n_0}}{|z-x|^{N+2\alpha}}dz\\&\le &0, \quad\mbox{in}\quad B_r(x_0)\cap \mathcal{A}_{k}.
\end{eqnarray*}
In Case b) we have use \equ{gamma1} and in Case c) we have use \equ{gamma2}, to get similar conclusion.
Then, by the Maximum Principle, we get
$$
\bar w_{n_0}(x_0)\leq \max\{ M_{n_0}, \sup_{x\in B_r(x_0)\cap \mathcal{A}_{k}^c}(w_{k_0}+w)\}.
$$
Using the same argument as in \equ{mpw}, we conclude that 
$$
\sup_{x\in B_r(x_0)\cap \mathcal{A}_{k}^c}(w_{k_0}+w)>
M_{n_0}
$$
does not hold and  therefore
 \begin{equation}\label{www}
 \bar w_{n_0}(x_0)=w_k(x_0)+w(x_0)\leq  M_{n_0}.
 \end{equation}
 Next, by the definition of $M_n$, we choose $x_1 \in \Omega\setminus A_{r/n_0}$ such that $w_k (x_1)=M_{n_0}$.
But then we have
$$w_k(x_0)+w(x_0)\geq w(x_0)=\frac{K_1}{2}V(0)r^\tau\quad \mbox{in Case a) and b)}$$
$$
\quad\mbox{and }\quad w_k(x_0)+w(x_0)\geq w(x_0)=\frac{K_1}{2}V(0)r^{\gamma+2\alpha} \quad \mbox{in Case c)}.$$
Thus,  by the asymptotic behavior of $v$, \equ{pp1} in Case a), \equ{gamma1} in Case b) and  \equ{gamma2} in Case c), we have
$$r^\tau\geq n_0^\tau C v(x_1)\quad\mbox{and}\quad r^{\gamma+2\alpha}\geq r^{\gamma/p} \geq n_0^{\gamma/p} C v(x_1)\quad\mbox{in Case c)}.$$
We recall that in Case a) $K_1=C(k^p-k)$, so from \equ{www}
\begin{equation}\label{4.2.4}
u(x_1)>(1+c_0)kv(x_1),
\end{equation} 
where $c_0>0$ is a constant, not depending on $x_0$ and increasing in $k$.
Now we repeat this process above initiating by $x_1$ and $k_1=k(1+c_0)$. 
 Proceeding  inductively,  we can
find a sequence $\{x_m\}\subset \mathcal{A}$ such that
$$u(x_m)>(1+c_0)^mkv(x_m),$$ which contradicts the common asymptotic behavior of $u$  and $v$.

In the Case b) and c) recall that $K_1=C(k-1)$  and, as before, we can proceed  inductively to
find a sequence $\{x_m\}\subset \mathcal{A}$ such that
$$u(x_m)>(k+m c_0)v(x_m),$$
which again contradicts the common asymptotic behavior of $u$  and $v$.
 \hfill $\Box$\\

\setcounter{equation}{0}
\section{Proof of our non-existence results}

In this section we prove our non-existence results. Our arguments are based on the construction of some special super and sub-solutions and some ideas used in Section \S 5. The main portion of our proof is based on the following proposition that we state and prove next.
\begin{proposition}\label{th 4}
Assume that $\Omega$ is an open, bounded and connected domain of
class $C^2$, $\alpha\in(0,1)$, $p>1$ and $f$ is nonnegative.
 Suppose that $U$ is a sub or  super-solution of (\ref{eteo1}) satisfying $U=0$ in $\Omega^c$ and  
\equ{1.7} for some $\tau\in(-1,0)$. Moreover, if $\tau>-\frac{2\alpha}{p-1}$,
assume there are numbers $\epsilon>0$ and $\delta>0$ such that,  in
case  $U$ is a sub-solution of (\ref{eteo1}), 
\begin{equation}\label{eq.2.11}
(-\Delta)^\alpha U(x)\le-\epsilon d(x)^{\tau-2\alpha} \quad
\mbox{or}\quad f(x)\ge\epsilon d(x)^{\tau-2\alpha}, \quad  \mbox{for } x\in
A_{\delta},
\end{equation}
and in  case  $U$ is a super-solution of (\ref{eteo1}),
\begin{equation}\label{eq.2.111}
(-\Delta)^\alpha U(x)\ge \epsilon d(x)^{\tau-2\alpha}\ \ \mbox{and}\ \
f(x)\le\frac\epsilon2d(x)^{\tau-2\alpha},\quad\mbox{for }  x\in A_{\delta}.
\end{equation}
Then there is no solution $u$ of (\ref{eteo1}) such that, in  case $U$ is a sub-solution,
\begin{eqnarray}\label{4.1}\nonumber
0<\liminf_{x\in\Omega,\
x\to\partial\Omega}u(x)d(x)^{-\tau}&\le&\limsup_{x\in\Omega,\
x\to\partial\Omega}u(x)d(x)^{-\tau}\\
&<&\liminf_{x\in\Omega,\
x\to\partial\Omega}U(x)d(x)^{-\tau}
\end{eqnarray}
or in  case  $U$ is a super-solution,
\begin{eqnarray}\label{4.0.1}\nonumber
0<\limsup_{x\in\Omega,\
x\to\partial\Omega}U(x)d(x)^{-\tau}&<&\liminf_{x\in\Omega,\
x\to\partial\Omega}u(x)d(x)^{-\tau}\\
&\le&\limsup_{x\in\Omega,\
x\to\partial\Omega}u(x)d(x)^{-\tau}<\infty.
\end{eqnarray}
\end{proposition}

We prove  this proposition by a contradiction argument, so we assume that   $u$ is a solution of (\ref{eteo1}) satisfying
(\ref{4.1}) or (\ref{4.0.1}), depending on the fact that $U$ is a sub-solution or a super-solution. Since  $f$ is non-negative we have that 
$u>0$ in $\Omega$ and by our assumptions on $U$, there is a constant $C_0\ge1$ so that, 
in case $U$ is a sub-solution
\begin{equation}\label{4.3}
C_0^{-1}\le u(x)d(x)^{-\tau}<U(x)d(x)^{-\tau}\leq C_0,\ \ x\in
A_{\delta}\ \ 
\end{equation}
and, in case $U$ is a super-solution
\begin{equation}\label{4.3s}
C_0^{-1}\leq U(x)d(x)^{-\tau}<u(x)d(x)^{-\tau}\le C_0,\ \ x\in
A_{\delta}.
\end{equation}
Here $\delta$ is decreased if necessary so that 
(\ref{eq.2.11}), (\ref{eq.2.111}),   (\ref{4.3}) and \equ{4.3s}  hold.
We define  \begin{equation}\label{pi} \pi_{k}(x)= \left\{
\arraycolsep=1pt
\begin{array}{lll}
 U(x)-k u(x),\ \ \ \ &
\mbox{in\ case } U \mbox{ is a sub-solution},\\[2mm]
u(x)-k U(x),& \mbox{in\ case } U \mbox{ is a super-solution},
\end{array}
\right.
\end{equation}
where $k\ge0$. In order to prove Proposition \ref{th 4}, we need the following two preliminary lemmas. 
\begin{lemma}\label{lm 3.1+}
Under the hypotheses of Proposition \ref{th 4}. If
$\mathcal{A}_{k}=\{x\in  \Omega\,/\,\pi_k(x)> 0\}\not=\O,$
for $k>1$. Then,
\begin{equation}\label{4.1.5+}
\partial \mathcal{A}_{k}\cap\partial \Omega\not=\O.
\end{equation}
\end{lemma}
The proof of this lemma follows the same arguments as the proof of 
 Lemma  \ref{lm 3.1} so we omit it.
\begin{lemma}\label{lm 3.2***} Under the hypotheses of Proposition \ref{th 4}. If
$\mathcal{A}_{k}\not=\O,$
 for $k>1$, then
 \begin{equation}\label{64.1.6}
 \sup_{x\in\Omega}\pi_k(x)=+\infty.
 \end{equation}
\end{lemma}
\noindent
{\bf Proof.} 
If (\ref{64.1.6}) fails, then we have $
M=\sup_{x\in \Omega}\pi_k(x)<+\infty.$ We see that $ M>0$ and, as in Lemma \ref{lm 3.2}, there is no point  $\bar x\in\Omega$ achieving $M$.  By Lemma \ref{lm 3.1+} we may choose $x_0\in\mathcal{A}_{k}$
and $r=d(x_0)/4$ such that $B_r(x_0)\subset A_\delta$, where $r$ could be chosen as small as we want. Here 
$\delta$ is as in (\ref{eq.2.11}) and (\ref{eq.2.111}). 

In what follows we consider  $x\in  B_r(x_0)\cap\mathcal{A}_{k}$ and we notice that  $3r<d(x)<5r$. We first analyze the case $U$ is a   sub-solution 
and $\tau\le-\frac{2\alpha}{p-1}$. We have
\begin{eqnarray*}
 (-\Delta)^\alpha
 \pi_{k}(x)&\leq&-U^p(x)+ku^p(x)-(k-1)f(x)
 \\&\le&-(k^{p-1}-1)k  u^p(x)
 \\&\le&-C^{-p}_0(k^{p-1}-1)k d(x)^{\tau p}
 \le-K_1r^{\tau-2\alpha},
 \end{eqnarray*}
where we have used  $f\ge0$, $k>1$,  (\ref{4.3}), $K_1=5^{\tau-2\alpha}C^{-p}_0(k ^{p-1}-1)k
>0$ and   $C_0$ is taken
from (\ref{4.3}).
Next we consider the case $U$ is a   sub-solution and $\tau>-\frac{2\alpha}{p-1}$. By the first inequality in
(\ref{eq.2.11}), we have
\begin{eqnarray*}
 (-\Delta)^\alpha
 \pi_{k }(x)&\leq&-\epsilon d(x)^{\tau-2\alpha}+k u^p(x)-kf(x)
 \\&\le&-(\epsilon-kC^{p}_0r^{2\alpha-\tau+\tau p}) d(x)^{\tau-2\alpha}\le -K_1 r^{\tau-2\alpha},
  \end{eqnarray*}
where the last inequality is achieved by choosing  $r$ small enough so that  $(\epsilon-kC^{p}_0r^{2\alpha-\tau+\tau
p})\ge\frac\epsilon2$ and $K_1=5^{\tau-2\alpha}\frac\epsilon2$.
On the other hand, if the second inequality in
(\ref{eq.2.11}) holds, we have
\begin{eqnarray*}
 (-\Delta)^\alpha
 \pi_{k }(x)&\leq&k u^p(x)-(k-1)\epsilon d(x)^{\tau-2\alpha}
 \\&\le&-((k-1)\epsilon-kC^{p}_0r^{2\alpha-\tau+\tau p}) d(x)^{\tau-2\alpha}\le -K_1 r^{\tau-2\alpha}, 
 \end{eqnarray*}
where $r$ satisfies 
$(k-1)\epsilon-kC^{p}_0r^{2\alpha-\tau+\tau
p}\ge\frac{k-1}2\epsilon$ and $K_1=5^{\tau-2\alpha}\frac{k-1}2\epsilon$.

In  case $U$ is a super-solution and 
$\tau\le-\frac{2\alpha}{p-1}$, we argue similarly to obtain
\begin{eqnarray*}
 (-\Delta)^\alpha
 \pi_{k }(x)&\leq&-u^p(x)+k U_1^p(x)-(k-1)f(x)
 \le-K_1r^{\tau-2\alpha},
 \end{eqnarray*}
where $K_1=5^{\tau-2\alpha}C^{-p}_0(k ^{p-1}-1)k >0$.
Finally, in case $U$ is a super-solution and  $\tau>-\frac{2\alpha}{p-1}$, using (\ref{eq.2.111}) we find
\begin{eqnarray*}
 (-\Delta)^\alpha
 \pi_{k }(x)&\leq&-u^p(x)-k\epsilon d(x)^{\tau-2\alpha}+ f( x)
 \le -K_1r^{\tau-2\alpha},
 \end{eqnarray*}
with $K_1=5^{\tau-2\alpha}\frac k2\epsilon >0$. Thus, in all cases we have obtained 
\begin{equation}\label{3.0.12+}
(-\Delta)^\alpha
 \pi_{k }(x)\leq-K_1 r^{\tau-2\alpha},\ \ x\in  B_r(x_0)\cap\mathcal{A}_{k},
\end{equation}
for some $K_1=K_1(k)>0$ non-decreasing with $k$.
From here we can argue as in Lemma \ref{lm 3.2} to get a contradiction 
$\Box$

 Now proof of Proposition \ref{th 4} is easy. 
 
\noindent{\bf Proof of Proposition \ref{th 4}.} From \equ{3.0.12+}, recalling that   $K_1$ non-decreasing with $k$, we can argue as in the proof of uniqueness result in Case b) to get
a sequence $(x_m)$ in $A_\delta$ such that, for some $k_0>1$ and $\bar k>0$,  in case $U$ is a sub-solution we have
 $$
U(x_m)>(k_0+m\bar k) u(x_m)\ 
$$
and, in case $U$ is a super-solution we have
$$
u(x_m)>(k_0+m\bar k) U(x_m).
$$
From here  we obtain a contradiction with (\ref{4.3}) or (\ref{4.3s}), for   $m$ large. 
\hfill $\Box$\\

\noindent{\bf Proof of  non-existence part of Theorem \ref{th 3.1.1}.}   For any $t>0$ we
 construct a sub-solution or  super-solution $U$ of
(\ref{eteo1}) such that
\begin{equation}\label{4.0.001}
\lim_{x\in\Omega,x\to\partial\Omega}U(x)d(x)^{-\tau}=t,
\end{equation}
and $U$ satisfies the assumption of Proposition  \ref{th 4}, for different combinations of the parameters $p$ and $\tau$.
For $t>0$ and $\mu\in\R$ we define   
\begin{equation}\label{eq 2.11}
U_{\mu,t}=tV_\tau+\mu  V_0\ \ \mbox{in}\ \ \R^N,
\end{equation}
where $V_0=\chi_\Omega$ is the characteristic function of $\Omega$ and $V_\tau$ is defined in \equ{3.3.1}. 
It is obvious that
(\ref{4.0.001}) holds for $U_{\mu,t}$
 for any
$\mu\in\R$.   To complete proof we show that  for any $t>0$, there is
$\mu(t)$ such that $U_{\mu(t),t}$ is a sub-solution or super-solution of (\ref{eteo1}), depending on the zone to which $(p,\tau)$ belongs.

\textbf{Zone 1:} We consider $p>1$ and $\tau\in (\tau_0(\alpha),0)$. By Proposition 
\ref{lm 3.3.1} $(ii)$, there exist
 $\delta_1>0$ and $C_1>0$ such that
\begin{equation}\label{4.4}
(-\Delta)^{\alpha}V_\tau(x)>C_1d(x)^{\tau-2\alpha},\ \ x\in
A_{\delta_1}.
\end{equation}
Combining with $(H2^*)$, for any $\mu>0$, there exists
$\delta_1>0$ depending on $t$ such that
$$(-\Delta)^{\alpha}U_{\mu,t}(x)+U_{\mu,t}^p(x)-f(x)>C_1td(x)^{\tau-2\alpha}-Cd(x)^{-2\alpha}\ge0,\ \ x\in
A_{\delta_1}.$$
On the other hand,
 since $V_\tau$ is of class $C^2$, $f$ is continuous in $\Omega$ and $\Omega\setminus A_{\delta_1}$ is compact,  there
exists $C_2>0$ such that
\begin{equation}\label{44.5}
|f|,\ |(-\Delta)^\alpha V_\tau(x)|\leq C_2,\ \ x\in\Omega\setminus
A_{\delta_1}.
\end{equation}
Then, using (\ref{V_0}),  there exists $\mu>0$ such that
\begin{equation}\label{4.5}(-\Delta)^{\alpha}U_{\mu,t}(x)+U_{\mu,t}^p(x)-f(x)>-2C_2+C_0\mu\ge0,\ \ x\in
\Omega\setminus A_{\delta_1}.
\end{equation} 
We conclude that for any
$t>0$, there exists $\mu(t)>0$ such that  $U_{\mu(t),t}$ is a super-solution of (\ref{eteo1}) and, by
$(H2^*)$ and \equ{4.4}, it satisfies (\ref{eq.2.111}).

\textbf{Zone 2:} We consider   $p> 1+2\alpha$ and $\tau\in(-1,-\frac{2\alpha}{p-1})$.  By Proposition \ref{lm 3.3.1} $(i)$ and $(ii)$, there exists
$\delta_1>0$ depending on $t$ such that
\begin{equation}
\label{4.700}
 (-\Delta)^\alpha
 U_{\mu,t}(x)+
 U^p_{\mu,t}(x)-f(x)\geq-C_1t d(x)^{\tau-2\alpha}+t^pd(x)^{\tau p}-Cd(x)^{-2\alpha}\ge0,
 \end{equation}
 for $x\in A_{\delta_1}$ and for any $\mu>0$, where we used that  $0>\tau-2\alpha>\tau
p$.  On the other hand, for $x\in \Omega\setminus A_{\delta_1}$,
(\ref{4.5}) holds for some $\mu>0$ and so we have constructed a super-solution of (\ref{eteo1}).

 \textbf{Zone 3:}  We consider $1+2\alpha<p\le1-\frac{2\alpha}{\tau_0(\alpha)}$ and
 $\tau\in(-\frac{2\alpha}{p-1},\tau_0(\alpha))$, which implies that $\tau
 p>\tau-2\alpha$.
  By Proposition \ref{lm 3.3.1} $(i)$ and $f\ge0$ in $\Omega$, there exists $\delta_1>0$ so that  for all
 $\mu\le0$
 \begin{equation}\label{eq 2.13}
 (-\Delta)^\alpha
 U_{\mu,t}(x)+U_{\mu,t}^p(x)-f(x)\leq-C_1t d(x)^{\tau-2\alpha}+t^pd(x)^{\tau
 p}\le0,
\end{equation}
 for $x\in
A_{\delta_1}$.
Then, using (\ref{V_0}) and (\ref{44.5}),  there exists
$\mu=\mu(t)<0$ such that
\begin{equation}\label{4.500}
(-\Delta)^{\alpha}U_{\mu,t}(x)+U_{\mu,t}^p(x)-f(x)<2C_2+C_0\mu\le0,\ \ x\in
\Omega\setminus A_{\delta_1}.
\end{equation} 
 We conclude that for any $t>0$, there
exists $\mu(t)<0$ such that  $U_{\mu(t),t}$  is a sub-solution of (\ref{eteo1}) and it satisfies
(\ref{eq.2.11}).

We see that Zone 1, 2 and 3 cover the range of parameters in part $(i)$ of Theorem \ref{th 3.1.1}, completing the proof in the case. 

\textbf{Zone 4:} To cover part (ii) of Theorem \ref{th 3.1.1} we only need to consider $p=
1-\frac{2\alpha}{\tau_0(\alpha)}$ with
  $\tau=\tau_0(\alpha)=-\frac{2\alpha}{p-1}$, which implies that  $\tau
  p=\tau-2\alpha<\min\{\tau-2\alpha+\tau+1,\tau\}$.
  By Proposition
\ref{lm 3.3.1}  $(iii)$, there exists $\delta_1>0$ depending on $t$
such that
\begin{eqnarray*}
 (-\Delta)^\alpha
 U_{\mu,t}(x)+
 U^p_{\mu,t}(x)-f(x)&\geq&-C_1t d(x)^{\min\{\tau-2\alpha+\tau+1,\tau\}}+t^pd(x)^{\tau p}\\& &-Cd(x)^{-2\alpha} \ge0,\ \ \ x\in A_{\delta_1}
\end{eqnarray*}
for any  $\mu>0$. For $x\in \Omega\setminus A_{\delta_1}$,
(\ref{4.5}) holds for some $\mu>0$, so  we have constructed a super-solution of (\ref{eteo1}).

We see that Zones 1, 2 and 4 cover the parameters in  part $(ii)$ of Theorem \ref{th 3.1.1}, so the proof is complete in this case too.  

 \textbf{Zone 5:}  We consider $1<p\le1+2\alpha$ and
 $\tau\in(-1,\tau_0(\alpha))$, which implies that $\tau
 p>\tau-2\alpha$.
  By Proposition \ref{lm 3.3.1} $(i)$ and $f\ge0$ in $\Omega$, there exists $\delta_1>0$ such that  for all
 $\mu\le0$ and  $x\in
A_{\delta_1}$, inequality \equ{eq 2.13} holds.
Then, using (\ref{V_0}) and (\ref{44.5}),  there exists
$\mu=\mu(t)<0$ such that \equ{4.500} holds and 
we conclude that for any
$t>0$, there exists $\mu(t)<0$ such that  $U_{\mu(t),t}$ satisfies
the first inequality of (\ref{eq.2.11}) and it is a sub-solution
of (\ref{eteo1}).

We see that Zones 1 and 5 cover the parameters in part $(iii)$ of Theorem \ref{th 3.1.1}. This completes the proof.\hfill$\Box$\\

\noindent{\bf Proof of  Theorem \ref{nonexistence}.} Here again we construct sub or super-solutions satisfying Proposition \ref{th 4} to prove the theorem.  In the case of a weak source, that is,
part $(i)$ of  Theorem \ref{nonexistence}, we have  $p\ge
1-\frac{2\alpha}{\tau_0(\alpha)}$ and
$-2\alpha-\frac{2\alpha}{p-1}\le \gamma<-2\alpha$, which implies
that $-1<\tau_0(\alpha)\le -\frac{2\alpha}{p-1}\le\gamma+2\alpha<0$. We consider two zones depending on $\tau$.

\textbf{Zone 1:} we consider $\tau\in (\gamma+2\alpha,0)$, so we have  $\gamma<\tau p$ and $\gamma<\tau -2\alpha$. By Proposition
 \ref{lm
3.3.1} $(ii)$ and $(H3)$, we have that,  for any $t>0$  there exist
 $\delta_1>0$, $C_1>0$ and $C_2>0$
 such that
\begin{equation}\label{eq 2.14}(-\Delta)^\alpha U_{\mu,t}(x)+U_{\mu,t}^p(x)-f(x) \le
C_1td(x)^{\tau-2\alpha}+t^pd(x)^{\tau p}-C_2d(x)^{\gamma}\le0,
\end{equation} 
for $ x\in A_{\delta_1}$ and  any $\mu\le0$. On the other hand, 
 using (\ref{V_0}) and (\ref{44.5}) we find
$\mu=\mu(t)<0$ such that (\ref{4.500}) holds for $ x\in
\Omega\setminus A_{\delta_1}.$  We conclude that for any $t>0$,
there exists $\mu(t)<0$ such that  $U_{\mu(t),t}$ is is a sub-solution of
(\ref{eteo1}) and by $(H3)$,
it satisfies (\ref{eq.2.11}).

\textbf{Zone 2:}  we consider $\tau\in(-1,\gamma+2\alpha)$. For
 $\tau\in(\tau_0(\alpha),\gamma+2\alpha)$ in case  $\tau_0(\alpha)<\gamma+2\alpha$,
by Proposition \ref{lm 3.3.1} $(i)$ there exists $\delta_1>0$, depending on $t$, such that 
\begin{equation}\label{eq4.700}
 (-\Delta)^\alpha
 U_{\mu,t}(x)+
 U^p_{\mu,t}(x)-f(x)\geq C_1t d(x)^{\tau-2\alpha}-C_2d(x)^{\gamma}\ge0,
\end{equation} 
for $x\in
 A_{\delta_1}$ and  any $\mu\ge0$. For
 $\tau\in(-1,\tau_0(\alpha)]\cap(-1,\gamma+2\alpha)$, we have
 $\tau p<\gamma$ and $\tau p<\tau-2\alpha$, so by Proposition \ref{lm 3.3.1} $(ii)$ and $(iii)$,
there exists $\delta_1>0$ dependent of $t$ such that (\ref{4.700})
holds for any $\mu\ge0$, while for $x\in \Omega\setminus A_{\delta_1}$, (\ref{4.5}) holds for some
$\mu>0$. We conclude that for any $t>0$, there exists $\mu(t)>0$
such that  $U_{\mu(t),t}$ is
a super-solution of (\ref{eteo1}) and  by $(H3)$ it satisfies (\ref{eq.2.111}), completing the proof in the weak source case.

Next we consider the case of strong source, that is part $(ii)$ of  Theorem \ref{nonexistence}. Here we
have that $$-1<\frac\gamma p<-\frac{2\alpha}{p-1}<0.$$
Here again we have two zones, depending on the parameter $\tau$.

\textbf{Zone 1:}
we consider $\tau\in(\frac\gamma p,0)$, in which case we have 
$\tau-2\alpha>\gamma$ and $\tau p>\gamma $. Then there exist
 $\delta_1>0$, $C_1>0$ and $C_2>0$
 such that (\ref{eq 2.14}) holds
 for any $\mu\le0$ and 
using (\ref{V_0}) and (\ref{44.5}),  there exists
$\mu=\mu(t)<0$ such that (\ref{4.500}) holds for $ x\in
\Omega\setminus A_{\delta_1}.$ Thus, for any $t>0$ there exists
$\mu(t)<0$ such that  $U_{\mu(t),t}$ is a  sub-solution of
(\ref{eteo1}) and $(H3)$  implies the first
inequality of  (\ref{eq.2.11}).  

\textbf{Zone 2:} we consider $\tau\in(-1,\frac\gamma p)$, in which case we have $\tau
p<\tau-2\alpha$ and $\tau p<\gamma $. Then there exist
 $\delta_1>0$, $C_1>0$ and $C_2>0$
 such that (\ref{eq4.700}) holds
 for $x\in A_{\delta_1}$ and $\mu\ge0$. We see also that for $x\in \Omega\setminus A_{\delta_1}$, inequality (\ref{4.5}) holds for some
$\mu>0$and so for any $t>0$, there exists $\mu(t)>0$ such that
$U_{\mu(t),t}$ is a super-solution of (\ref{eteo1}).

This completes the proof of the theorem.\hfill$\Box$


\begin{thebibliography}{99}

\bibitem {AR}
J. M. Arrieta and A. Rodr\'{i}guez-Bernal, Localization on the
boundary of blow-up for reaction-diffusion equations with nonlinear
boundary conditions, {\it Comm. Partial Diff. Eqns.}, 29, 1127-1148,
2004.


\bibitem {BM1}
C. Bandle and M. Marcus, Large solutions of semilinear elliptic
equations: Existence, uniqueness and asymptotic behaviour, {\it J.
Anal. Math.}, 58, 9-24, 1992.

\bibitem {BM2} C. Bandle and M. Marcus, Dependence of blowup rate of large
solutions of semilinear elliptic equations on the curvature of the
boundary, {\it Complex Variables Theory Appl.}, 49, 555-570, 2004.

\bibitem {CC} X. Cabr\'{e} and L. Caffarelli, Fully Nonlinear Elliptic Equation,
{\it American Mathematical Society, Colloquium Publication,} Vol.
43, 1995.

\bibitem {CS0} L. Caffarelli, S. Salsa and L. Silvestre,
Regularity estimates for the solution and the free boundary to the
obstacle problem for the fractional Laplacian, {\it Inventiones
mathematicae}, 171, 425-461, 2008.

\bibitem {CS1} L. Caffarelli and L. Silvestre, Regularity theory for fully nonlinear integro-differential equaitons,  {\it
Comm. Pure Appl. Math.}, 62(5), 597-638, 2009.

\bibitem {CS2} L. Caffarelli and L. Silvestre, Regularity results for nonlocal
equations by approximation, {\it Arch. Ration. Mech. Anal.}, 200(1),
59-88, 2011.

\bibitem {CS_Evans} L. Caffarelli and L. Silvestre,  The Evans-Krylov theorem for non local fully non linear equations.
Annals of Mathematics.174 (2), 1163-1187, 2011. 

\bibitem {CF2} H. Chen and P. Felmer, Liouville Property for fully nonlinear integral equation
in exterior domain, {\it Preprint.}

\bibitem {CCE}
M. Chuaqui, C. Cort\'{a}zar, M. Elgueta and J.
Garc\'{i}a-Meli\'{a}n, Uniqueness and boundary behaviour of large
solutions to elliptic problems with singular weights, {\it Comm.
Pure Appl. Anal.}, 3, 653-662, 2004.

\bibitem {DL} M. del Pino and R. Letelier, The influence of domain geometry in
boundary blow-up elliptic problems, {\it Nonlinear Analysis: Theory,
Methods \& Applications}, 48(6), 897-904, 2002.

\bibitem {DL1} G. D\'{i}az and R. Letelier, Explosive solutions of quasilinear
elliptic equations: existence and uniqueness, {\it Nonlinear
Analysis: Theory, Methods \& Applications,} 20(2), 97-125, 1993.

\bibitem {DH}
Y. Du and Q. Huang, Blow-up solutions for a class of semilinear
elliptic and parabolic equations, {\it SIAM J. Math. Anal.}, 31,
1-18, 1999.


\bibitem {FQ0} P.  Felmer and A. Quaas, Boundary blow up solutions for fractional elliptic equations. {\it Asymptotic Analysis}, Volume 78 (3), 123-144,
 2012.
 
 \bibitem {FQT} P. Felmer, A. Quaas, J. Tan, Positive solutions of nonlinear Schrodinger equation with the fractional Laplacian, Proceedings of the Royal Society of Edinburgh: Section A Mathematics, 142, 1-26, 2012.


\bibitem {FQ1} P.  Felmer and A. Quaas, Fundamental solutions and
two properties of elliptic maximal and minimal operators, {\it
Trans. Amer. Math. Soc.}, 361(11), 5721-5736, 2009.

\bibitem {FQ2} P. Felmer and A. Quaas, Fundamental solutions and Liouville type theorems for nonlinear integral operators, {\it
Advances in Mathematics}, 226, 2712-2738, 2011.

\bibitem {FQ3} P. Felmer and A. Quaas,  Fundamental solutions for a class of Isaacs integral operators, {\it
Discrete and Continuous Dynamical Systems}, 30(2), 493-508, 2011.

\bibitem {G}  J. Garc\'{i}a-Meli\'{a}n, Nondegeneracy and uniqueness for boundary blow-up elliptic
problems, {\it J. Diff. Eqns.}, 223(1), 208-227, 2006.


\bibitem {GG} J. Garc\'{i}a-Meli\'{i}an, R. G\'omez-Re\~nasco, J.
L\'opez-G\'omez and J. Sabina de Lis, Pointwise growth and
uniqueness of positive solutions for a class of sublinear elliptic
problems where bifurcation from infity occurs, {\it Arch. Ration.
Mech. Anal.}, 145(3), 261-289, 1998.

\bibitem {I} H. Ishii, On uniqueness and existence of viscosity
solutions of fully nonlinear second-order elliptic PDE's, {\it Comm.
Pure Appl. Math.}, 42(1), 15-45, 1989.

\bibitem {K} J. B. Keller, On solutions of $\Delta u = f(u)$, {\it Comm. Pure Appl.
Math.},  10,  503-510, 1957.

\bibitem {K1} S. Kim, A note on boundary blow-up problem of $\Delta u = u^p$, {\it IMA preprint
No.}, 18-20, 2002.

\bibitem {LN} C. Loewner and L. Nirenberg, Partial differential equations invariant
under conformal projective transformations, in Contributions to
Analysis (a collection of papers dedicated to Lipman Bers). {\it
Academic Press, New York,} 245-272, 1974.

\bibitem {MV} M. Marcus and L. V\'{e}ron, Existence and uniqueness results for large solutions of
general nonlinear elliptic equation, {\it J. Evol. Equ.} 3, 637-652,
2003.

\bibitem {MV1} M. Marcus and L. V\'{e}ron, Uniqueness and asymptotic behavior of solutions
with boundary blow-up for a class of nonlinear elliptic equations,
{\it Ann. Inst. H. Poincar\'{e}} 14(2),  237-274, 1997.


\bibitem {O} R. Osserman, On the inequality $\Delta u = f(u)$, {\it Pac. J. Math.} 7,
1641-1647, 1957.


\bibitem{PVS} G. Palatucci, O. Savin and  E. Valdinoci, Local  and global
minimizers for a variational energy involving a fractional norm,
{\it http://arxiv.org/abs/1104.1725.}




\bibitem {R} V. R\v{a}dulescu, Singular phenomena in nonlinear elliptic
problems: from  blow-up boundary solutions to equations with
singular nonlinearities, {\it Handbook of Differential Equations:
Stationary Partial Differential Equations}, 4, 485-593, 2007.

\bibitem{rosoton} X. Ros-oton and J. Serra, The Dirichlet problem for the fractional laplacian, regularity up to the boundary, {\it http://arxiv.org/abs/1207.5985.}

\bibitem {S} L. Silvestre, H\"{o}lder estimates for solutions of integro differential
equations like the fractional laplace. {\it Indiana Univ. Math. J.},
55, 1155-1174, 2006.

\bibitem{STEIN} E.M. Stein, Singular Integrals and Differentiability Properties of Functions,
Princeton University Press, 1970.


\bibitem {V} L. V\'{e}ron, Semilinear elliptic equations with uniform blow-up on the
boundary, {\it J. Anal. Math.}, 59(1), 231-250, 1992.

\bibitem {Z} Z. Zhang, A remark on the existence of explosive solutions for a
class of semilinear elliptic equations, {\it  Nonlinear Analysis:
Theory, Methods \& Applications}, 41, 143-148, 2000.
\end{thebibliography}
\end{document}